\documentclass[journal]{IEEEtran}

\usepackage{array}
\usepackage{algorithm}
\usepackage[noend]{algpseudocode}
\usepackage{amsmath}
\usepackage{caption}
\usepackage{subcaption}

\usepackage{cases}
\usepackage{graphicx}
\usepackage{float}
\usepackage{mathtools}
\usepackage{ifpdf}
\usepackage{hyperref}

\hypersetup{pdflinkmargin=3pt}

\makeatletter
\newcommand*{\rom}[1]{\expandafter\@slowromancap\romannumeral #1@}
\makeatother

\begin{document}
%
\title{Vehicle Risk Assessment and Control for Lane-Keeping and Collision Avoidance in Urban and Highway Driving Scenarios}

\author{Hazem~M.~Fahmy,~\IEEEmembership{}
        Mohamed~A.~Abd El Ghany,~\IEEEmembership{IEEE Senior Member}
        and~Gerd~Baumann~\IEEEmembership{}
\noindent \thanks{ }
\thanks{Hazem M. Fahmy is with the Department of Electronics Engineering, German University in Cairo, Egypt.}
\thanks{Gerd Baumann is with the Department of Mathematics, German University in Cairo, Egypt. }
\thanks{Mohamed A. Abd El Ghany is with the Department of Electronics Engineering, German University in Cairo, Egypt and TU Darmstadt, Germany.}}

\markboth{}%
{Fahmy \MakeLowercase{\textit{et al.}}: Vehicle Risk Assessment and Control for Lane-Keeping and Collision Avoidance at Low-Speed and High-Speed Scenarios}

\maketitle

\begin{abstract}
This article examines a symbolic numerical approach to optimize a vehicle’s track for autonomous driving and collision avoidance. The new approach uses the classical cost function definition incorporating the essential aspects of the dynamic state of the vehicle as position, orientation, time sampling, and constraints on slip angles of tires. The optimization processes minimize the cost function and simultaneously determine the optimal track by varying steering and breaking amplitudes. The current velocity of the vehicle is limited to a maximal velocity, thus, allowing a stable search of the optimal track. The parametric definition of obstacles generates a flexible environment for low and high speed simulations. The minimal number of influential optimization variables guarantees a stable and direct generation of optimal results. By the current new approach to control a vehicle on an optimal track, we are able to autonomously move the vehicle on an arbitrary track approximated by low order polynomials. The optimization approach is also able to deal with a variety of different obstacles and the corresponding optimal smooth obstacle path. The computations demonstrate the effective control of a four wheel vehicle in normal operation and exceptional obstacle avoidance with continuously differentiable obstacle avoidance tracks. Simulation tests are done using vehicle’s velocities of 3m/s, 6m/s, 7.6m/s, 10m/s, 12 m/s, and 18m/s. At higher vehicle’s velocities, a mathematical-only approach is not sufficient and a mechanical intervention for tires is needed as a complimentary part to control the slip angle. The results shows that the cost function reached a considerably high average convergence-to-zero rate success in most of the tested scenarios.
\end{abstract}


\section*{Nomenclature}
\begin{description}[\setlabelwidth{$\alpha$}\usemathlabelsep]
\item[$\alpha$] \qquad  Tire slip angle [rad].
\item[$\beta_r$] \qquad  Braking ratio [-].
\item[$\delta$] \qquad  Steering angle [rad].
\item[$\dot{x}$] \qquad  Longitudinal velocity [m/s].
\item[$\dot{y}$] \qquad  Lateral velocity [m/s].
\item[$\dot{\psi}$] \qquad  Turning rate [rad/s].
\item[$F_x$] \qquad  Longitudinal tire force in the vehicle frame [N].
\item[$f_x$] \qquad  Longitudinal tire force in the tire frame [N].
\item[$F_y$] \qquad  Lateral tire force in the vehicle frame [N].
\item[$f_y$] \qquad  Lateral tire force in the tire frame [N].
\item[$m$] \qquad  Vehicle mass [m].
\item[$w_t$] \qquad  Vehicle track width [m].
\item[$l_f$] \qquad  Distance from vehicle center of gravity (CoG) to front axle [m].
\item[$l_r$] \qquad  Distance from vehicle CoG to rear axle [m].
\item[$J_z$] \qquad  Vehicle yaw inertia [kg.$m^2$].
\item[$C_\alpha$] \qquad  Tire cornering stiffness [-].
\item[$\mu$] \qquad  Road friction coefficient [-].
\item[$\epsilon$] \qquad  Road obstacle width [m].
\item[$F_z$] \qquad  Vertical load at each wheel [m].
\item[$N_d$] \qquad  Number of divisions for the time interval [-].
\item[$p_i$] \qquad  Penalty weight for each i component [-].
\item[$v_{max}$] \qquad  Maximum longitudinal velocity to be maintained

\end{description}

%
\IEEEpeerreviewmaketitle

\section{Introduction}
\IEEEPARstart{T}{he} European Road Safety Observatory estimated in their annual report \cite{euro} in 2015 that 26,000 fatalities are caused by road accidents each year within the European Union only. For every death on Europe's roads there are an estimated 4 permanently disabling injuries such as damage to the brain or spinal cord, 8 serious injuries and 50 minor injuries. In
2003, The National Highway Traffic Safety Administration (NHTSA) has reported in \cite{nhtsa} that among the various crash
types the most frequent one is the rear-end collision. Police reports in the United States supported the claim by showing that 29\% of the crash reports are caused by rear-end collisions. Previous study \cite{ger} showed that unintentional lane departure accidents were accounted for 14\% of all accidents reported in Germany only which also caused 30\% of all road accidents fatalities in the year 2013. Moreover, a research has predicted that road traffic accidents will increase globally by 67\% by the year 2020 \cite{research}. And recently in 2016, the World Health Organization (WHO) has stated a key fact that 1.25 million people die each year as a result of road traffic crashes, while road traffic injuries are the leading cause of death among people aged 15-29 years \cite{WHO}.

One way to reduce the number of accidents and their causalities is to actively assist road users in their driving. This is called preventive or active safety, and it ranges from Electronic Stability Control Systems to Driver Drowsiness-detection Systems.  Examples of advisory systems are dynamic active display, and driver drowsiness monitoring were discussed in \cite{3}, and \cite{8}. Other papers have discussed mechanical intervene  systems such as electronics stability systems \cite{11}, and anti-lock braking systems \cite{7}. Some studies also have presented threat assessment systems such as lane-keeping systems \cite{4,9}, and rear-end collision avoidance systems \cite{24}. In such systems, a preventive action is taken in a pre-accident phase either by giving alarms to the driver or assessing the situation and intervening in the vehicle mechanical behavior.

On one hand, Active Safety Systems are a complete frame work where a threat assessment, decision-making, and intervention modules are integrated into one block. The Active Safety System monitors the wiring of all these modules with each other. The threat assessment block is responsible for collecting the environment data by sensors fusion technique then applying mathematical modeling for estimating the risk probability \cite{24}. Decision making block is the mind of the process where optimization and control take place to decide whether an intervention is needed or not based on the data collected by the former block. An approach that studied such kind of algorithmic block was presented in \cite{47} as Neural Networks and in \cite{20} as Model Predictive Control. The final block, intervention module, is mainly a mechanical hardware related module where the intervention is done by manipulating throttle, tires or steering wheels. Vehicle Risk Assessment research in \cite{r1,r2,r3,r4,r5} discuss various aspects of estimating the risk for vehicle motion. On the other hand, the previous work for Passive Safety systems \cite{r6} had limited reliability in predicting the likelihood of having an accident as it only took action after the accident's occurrence. In active safety systems, some presented work in \cite{3} was only able to give the driver an advisory signal which will not contribute in preventing an accident in an effective way since the driver could ignore the advisory signal and an accident could still occur. Other active safety systems \cite{4} were designed to maintain the vehicle in the lane. These systems did not implement a system that can avoid an obstacle in such lane. The likelihood of having an accident will be high in such case which is not reliable in keeping the vehicle in a safe state. More active safety systems were designed to control the physical movement of the vehicle without including a prediction system or a specific algorithm. Such systems require a complimentary algorithm to decide the needed action to be taken.

\begin{figure}[t!]
     \centering
      \includegraphics[height = 0.3 \textheight , width=0.45\textwidth]{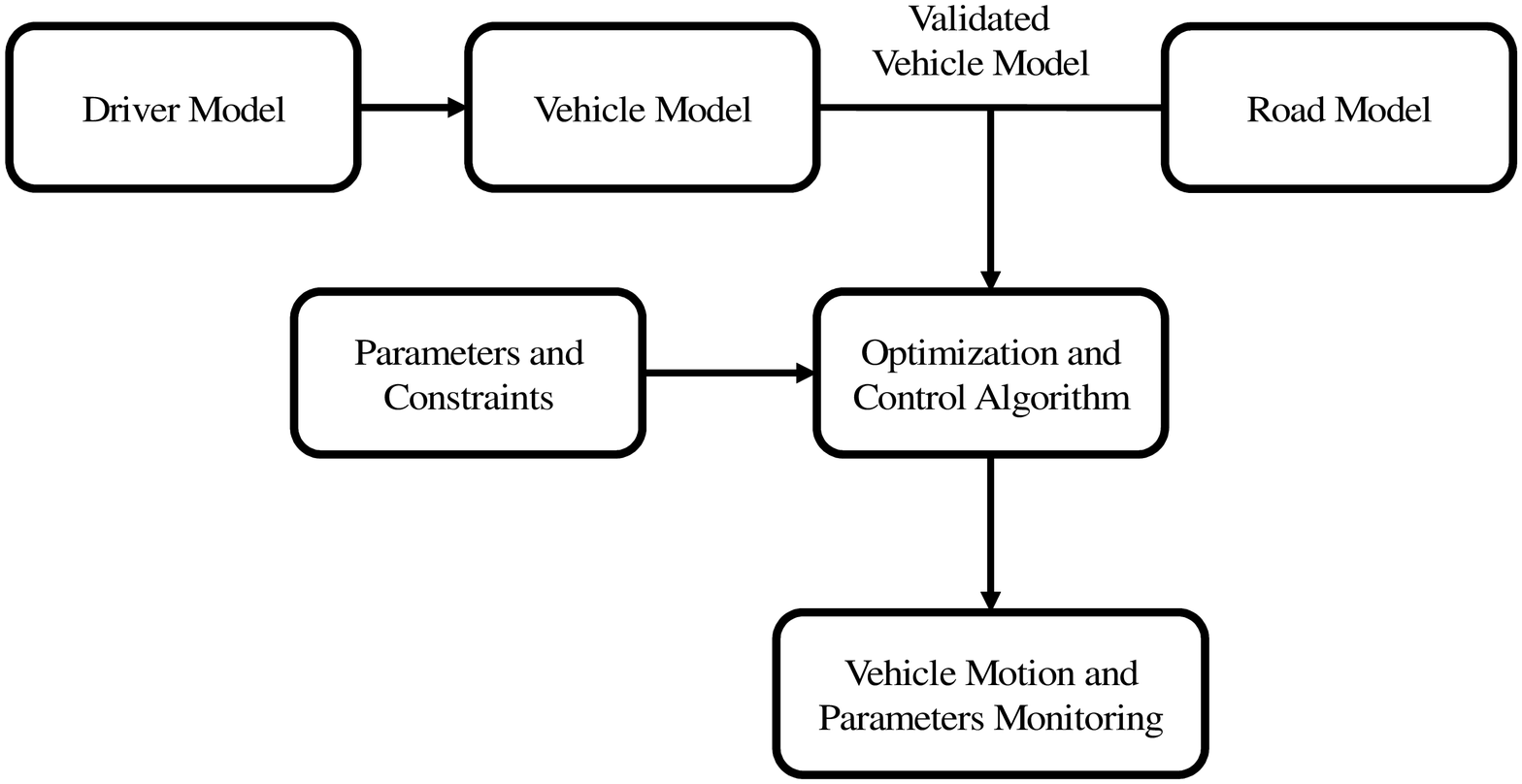}
    \caption{Block diagram illustrating the proposed framework integrity. The relation between vehicle model and driver model results into a vehicle path which is bound to the road model.}
		
      \label{a1}
\end{figure}



This paper is objected towards developing a risk-assessment algorithm that could control a vehicle to keep the presented lane and avoid a collision that may be caused by a road object. The algorithm is designed in a way that deals with the provided optimization problem in a complete symbolic numerical behavior. The constraints are set in order to meet a successful control behavior of the vehicle in terms of safety and stability. Furthermore, if a risk is identified, an action will continuously be taken until the error demolishes or tends to a minimal value. This action is done through controlling driving wheels and vehicle’s throttle. The whole approach can then be enclosed into an autonomous-driving system, that can automatically drive the vehicle into safe state based on the current environmental variables. The paper discusses a fully numerical approach through modeling, simulation, trial and error. The validation and simulation are done on a wide range of cases and it is included and discussed later in Section III. However, only interesting and useful results were included due to space limitations.

The use of numerical calculations through this work allows to extend the symbolic calculation if an analytical solution is
impossible. The combination of both numerical and symbolic approaches simplify first the formulation of the numerical methods and second allow to increase the effectiveness of a numerical algorithm. The mathematical modeling of the vehicle is needed in order to describe the motion according to the dynamics of the vehicle body, tires, and steering angle. The block diagram is shown in Fig. \ref{a1} that involves vehicle model and its validation through using driver model. The diagram also includes the road model, the optimization algorithm and its constraints. The monitoring block is used to finally see the resultant output of the optimization and control algorithm which is presented in Section \rom{3} discussing final results. A driver model is necessary for the validation of the vehicle behavior before going into the optimization and control phase. For the vehicle to be tested, a mathematical description of the road environment has to be introduced for the vehicle to react on such model. Test cases and scenarios are basically dependent on the road model geometry. In this study, the road can be modeled as a bent-road, straight-road or a road with an obstacle, and the optimization algorithm have to control the vehicle in following that pre-designed road model. All these blocks are then connected to the mind of the work or the optimization process. The optimization process is a classical approach using the dynamical behavior of the vehicle and driver model, or steering angle data, to achieve a minimization of a cost function $C$. The cost function is a representation of the amount of error in the process of achieving lane-keeping and collision avoidance. Some dynamic state variables are used during the optimization process for monitoring the status of the vehicle such as; slip angle, maximum velocity, braking ratio limits, steering angle limits. Such variables do not affect the cost function itself. However, the optimization algorithm is set to meet these constraints and not violate them while reaching a minimum cost function value. The function concludes all the parameters that can affect the lane-keeping or collision avoidance process achievement. The cost function is extended by requiring for specific optimization components a certain constraint. The constraints are inequality constraints limiting a subset of optimization variables to predefined intervals. The constraints and the use of penalty functions guarantee that a minimum can be found within the defined scope.

The rest of the paper is organized as follows; In Section \rom{2} the vehicle's dynamics and equations of motion that represent vehicle's behavior in addition to the road model are introduced. Section \rom{2}-C illustrates the mathematical behavior of the optimization and control algorithm in addition to the presented constraints and penalization technique. In Section \rom{3}, validation results of the proposed framework under low velocity and high velocity scenarios - including obstacles are - presented. In Section \rom{4} the results of the work is concluded and analyzed providing a summary of the whole proposed framework including remarks and outlining future work.

\section{Modeling}
In this section, the defining model equations used by the optimization and control algorithm is presented. Section \rom{2}-A introduces the used vehicle mathematical model in the work, Section \rom{2}-B introduces the modeled road geometry, and Section \rom{2}-C introduces the proposed optimization and control algorithm for lane keeping and collision avoidance combined.

\subsection{Vehicle Model}
For this study, a four-wheel vehicle model is adapted from the work \cite{20} which allows controlling and simulating the vehicle's dynamics due to the degrees-of-freedom of this model.

\begin{figure}[t!]
     \centering
      \includegraphics[scale = 0.4]{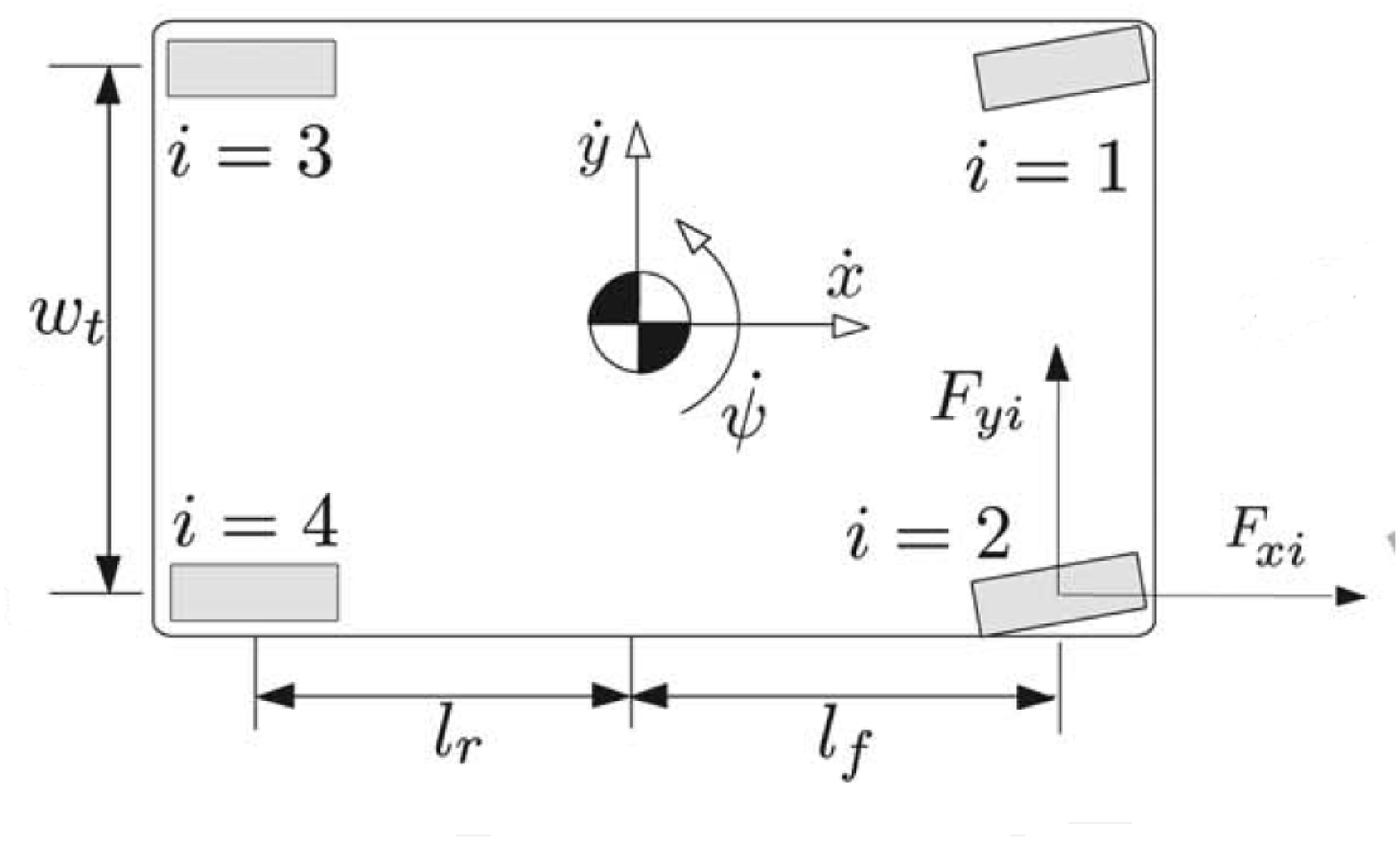}
    \caption{Modeling notation depicting the forces in the vehicle body fixed frame, the forces in the tire fixed frame, and the rotational and translational velocities \cite{20}.}
      \label{fig:1}
\end{figure}

The following set of $2^{nd}$ order differential equations for the variables $x$, $y$, $\psi$ are used to describe vehicle motion considering the vehicle sketch in Fig. \ref{fig:1}. Shown in Fig. \ref{fig:1} are the modeling notations that depict the forces in the vehicle body fixed frame, the forces in the tire fixed frame, and the rotational and translational velocities. Shown the illustration of \(\delta\) steering angle which is the angle between the vehicle orientation and the front wheels orientation. \(w_t , l_f,\) and \(l_r\) denote the dimension of the vehicle frame as the width, front, and rear lengths respectively. \(F_{xi}\) and \(F_{yi}\) denote the longitudinal and lateral forces respectively.

   \begin{equation}
m \ddot{x}  =  m \dot{y}\dot{\psi} + \sum\limits_{i=1}^4 F_{xi}
\end{equation}
\begin{equation}
m \ddot{y}  =  -m \dot{x}\dot{\psi} + \sum\limits_{i=1}^4 F_{yi}
\end{equation}
\begin{multline}
J_z \ddot{\psi}  =  l_f (F_{y1} + F_{y2}) - l_r (F_{y3} + F_{y4})+ \\
 \frac{w_t}{2} (-F_{x1} + F_{x2} - F_{x3} + F_{x4})
\end{multline}

\noindent where $x$, $y$ and $\psi$ are the longitudinal, lateral and turning positions of the vehicle respectively. $m$ is the vehicle's mass which is assumed to be constant throughout this work. $l_f$ and $l_r$ are the lengths of front and rear vehicle axles, $w_t$ is the width of the vehicle, $J_z$ is the vehicle yaw inertia, $F_{xi}$ and $F_{yi}$, for $i \in \{1,2,3,4\}$, are the longitudinal and lateral tire force components in the vehicle body frame for each tire and they are modeled as

\begin{equation}
F_{xi} = f_{xi} \cos(\delta_i) - f_{yi} \sin(\delta_i)
\end{equation}
\begin{equation}
F_{yi} = f_{xi} \sin(\delta_i) + f_{yi} \cos(\delta_i)
\end{equation}

\noindent where $\delta_i$  is the steering angle corresponding to wheel $i \in \{1,2,3,4\}$ and $f_{xi}$ is the longitudinal tire force in the tire frame which is computed in the following way:
    
    \begin{equation}
    f_{xi} = \beta_r \mu_i F_{zi} 
    \end{equation}
    
\noindent where $F_{zi}$ is the vertical tire force in vehicle frame for every tire - or the vertical load of every tire, and $\beta_r \in [-1,1]$ is referred to as the braking ratio. $\beta_r = -1$ corresponds to full braking while $\beta_r = 1$ corresponds to full throttle. Braking and throttle are taken into consideration as accelerating and decelerating so that the velocity will never be constant while $\beta_r \neq 0$.
    
 \noindent The lateral tire force in the tire frame is computed using a modified nonlinear Fiala tire model \cite{e} as shown here in (7a) and (7b) where 
    
    \begin{equation}
		f_{yi} = 
 \begin{cases} 
      -C_{\alpha_i} \tan(\alpha_i) + \frac{C_{\alpha_i}^2}{3 \eta \mu_i F_{zi}}
|\tan(\alpha_i)| \tan(\alpha_i) \\ - \frac{C_{\alpha_i}^3}{27 \eta^2 \mu_i^2
F_{zi}^2 } \tan(\alpha_i)^3 , \quad \alpha_i < \alpha_{sl} & (7a)\\
     -\eta \mu_i F_{zi} sgn(\alpha_i), \quad \quad \alpha_i \geq \alpha_{sl} & (7b)
   \end{cases}
\end{equation}

\noindent with $\alpha_i$, for $i \in \{1,2,3,4\}$, is the slip angle that can be modeled as an approximation as
    
    \begin{equation}
    \label{eq:alpha1}
    \alpha_1 = \alpha_2 = \frac{\dot{y} +l_f \dot{\psi}}{\dot{x}} - \delta
    \end{equation}
        \begin{equation}
        \label{eq:alpha3}
        \alpha_3 = \alpha_4 = \frac{\dot{y} -l_r \dot{\psi}}{\dot{x}}
    \end{equation}
     
\noindent The limits of the tires slip angle $\alpha$ were estimated as
		
		\begin{equation}
		 -\frac{2}{\pi} < \alpha_i < \frac{2}{\pi}
		\end{equation}
    
\noindent $\mu_i$ denotes the friction coefficient of the road which ranges from $[0,1]$ where $0$ corresponds to wet road while $1$ corresponds to a rough road, and  $F_z$ denotes the vertical load at each wheel. $C_{\alpha}$ is defined as the tire cornering stiffness, and $\alpha_{sl}$ is the slip limit angle and it is calculated as 
    
    \begin{equation}
\alpha_{sl} = \tan^{-1}(\frac{3 \mu \eta F_z }{C_\alpha})
\end{equation}
    
\noindent where $\eta$ is a braking dependent coefficient which is calculated as 
    
    \begin{equation}
\eta = \arctan(\frac{\sqrt[]{\mu^2 F_z^2 - f_x^2}}{\mu F_z})
\end{equation}
    
\noindent and it can be simplified to as
    
   \begin{equation}
\eta = \sqrt[]{1 - \beta_r^2}
\end{equation} 
   
    \textit{Assumption 1}: Since the steering angle at rear wheels are not controlled, $\delta_3, \delta_4$ are assumed to be zero, i.e. $\delta_3 = \delta_4 = 0$. Nevertheless, $\delta_{1,2}$ at front wheels can be controlled and have a value. It is assumed that $\delta_1 = \delta_2$, i.e., $\delta_1 = \delta_2 = \delta$.
    
    \textit{Assumption 2}: Vertical forces $F_{zi}$ are assumed constant and determined by the vehicle’s steady-state weight distribution when no lateral or longitudinal accelerations act at the vehicle center of gravity.
    
    \textit{Assumption 3}: The friction coefficient $\mu$ is assumed to be known and the same at all wheels, i.e., $\mu_i = \mu$, $\forall i$, and constant over a finite time horizon.
		
		\textit{Assumption 4}: In Eq. \ref{eq:alpha1} and \ref{eq:alpha3}, the vehicle’s longitudinal velocity $\dot{x}$ is assumed to never cross zero and never settles to a steady-state or the simulation is stopped due to a division by zero.

\subsection{Road Geometry Model}

In this section, road modeling is done using mathematical interpolation of road-line sections. Sectioning is done according
to the curvature of the road, so that an interpolation with a higher order polynomial function is performed if there is an object to avoid. The interpolated object-curve will be created using the dimensions of the obstacle being detected by sensors. $2^{nd}$ order polynomial function is followed if there is no object to avoid, where then the obstacle’s dimensions parameters will be set to zero. Listed below are the types of road models which were designed throughout the simulation process:

\subsubsection{Lane-Keeping Road Geometry} This geometry is described using a $2^{nd}$ order polynomial equation as follows,

\begin{equation}
R_{road}(x) = kx^2 + m
\end{equation}

\noindent where $k$ and $m$ are parameters that define the slope of the road in addition to the steepness and elevation, while $x$ is the longitudinal position of the vehicle. The geometry was formulated to be $x$-dependent instead of being time-dependent so that the road curvature is aligned with the $x$, $y$ coordinates graphing of the vehicle. A feature is added to the geometry in order to create a collision avoidance curve that will be discussed later in this section.

\begin{figure}[t!]
     \centering
      \includegraphics[scale = 0.8]{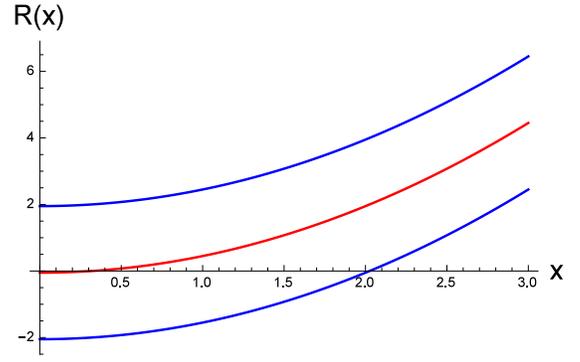}
    \caption{Vehicle lane simulation having an ideal centerline and two lane borders having no obstacles (lane-keeping).}
      \label{fig:2}
\end{figure}

\begin{figure}[t!]
     \centering
      \includegraphics[scale = 0.8]{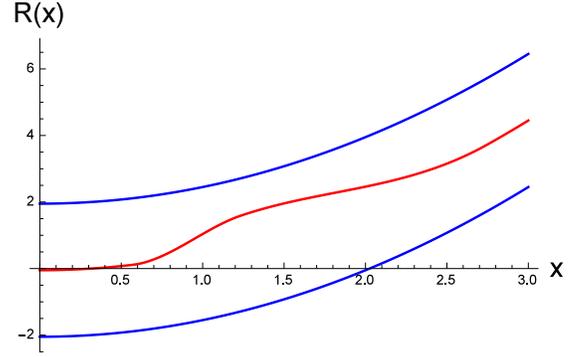}
    \caption{Vehicle lane simulation having an ideal centerline and two lane borders and an obstacle (obstacle-avoidance).}
      \label{fig:3}
\end{figure}

The resultant simulation of such road model is presented in Fig. \ref{fig:2} with  the ideally expected vehicle position in terms of longitudinal and lateral orientation $x$, and $y$.

\subsubsection{Collision-Avoidance Road Geometry} A collision-avoidance geometry can be described using the starting and ending points of the lane-keeping geometry in addition to an obstacle width parameter. The collision-avoidance geometry is treated as a double solution to the lane-keeping and collision-avoidance problems. Interpolation using Hermite technique was carried out at the sense of an obstacle that was parametrized using $\epsilon$ parameter where a constraint is applied such that:

\begin{equation}
\frac{-(x(t_{e}) - x_0)}{2} \leq \epsilon \leq \frac{x(t_{e}) - x_0}{2}
\end{equation}

\noindent where $x_0$ denotes the initial position of the obstacle, $x(t_{e})$ denotes the final position of the obstacle, and $\epsilon$ denotes the obstacle's width. The data of the obstacle can be extracted using sensor fusion technologies presented in \cite{b}. $\epsilon \neq 0$ is a required condition for the obstacle geometry to be created. The obstacle geometry is currently a rectangular box with width and depth, where width is denoted as final-initial coordinates, while depth is denoted as $\epsilon$. At the case where $\epsilon = 0$, the interpolation will fail and the road geometry would be a line-keeping path not a collision-avoidance path.

\noindent Interpolation was carried out using the information of the obstacle utilizing the interpolation technique of three-points describing the obstacle length and width and generating a slope for such curve. The resultant of collision-avoidance road model is plotted in Fig. \ref{fig:3} which is graphed for example in an interval $x \in [0,5]$, and $\epsilon = 0.5$. Noted that this is only an example for the purpose of illustrating obstacle geometry generation. $\epsilon$ can be positive or negative with regard to the position of ideal vehicle path as a reference zero. Positive $\epsilon$ denotes that the obstacle width would be upwards the center. While negative $\epsilon$ denotes that the obstacle width would downwards the center.

The flexible approach of 	\textit{Mathematica} was used to apply the most suitable technique for interpolation according to the nature of the curved road. The available techniques switches between the Spline and Hermite techniques in an automatic way.

\subsection{Optimization and Control Algorithm}
Optimization and control is needed for keeping the vehicle within the lane safe margins while having feedback information during the discrete time samples. The aim of optimization is to follow a known ideal line on the road. In addition to location and orientation some dynamic state variables like braking and throttle, slip angles and the maximum velocity are variables which can be used to characterize the current dynamical state.


The objective function, which is mentioned as the goal optimization process formulated in a mathematical way, is illustrated in this sub-section. The optimization parameters which control and affect the  objective function are discussed. Moreover, the constraints and penalty function that define the relation between the violation of the constraints and the objective function is also discussed in details in this sub-section.

\subsubsection{Objective Function}

The formulated problem has an objective function \(C\) to be met with a minimum cost while having constraints that should not be violated. The optimization problem has three varying field-of-search parameters, \((\delta ,\beta ,\text{$\delta $t})\), while five constraints to be fulfilled. Two of these parameters describe the vehicle's behavior; i.e. steering angle \(\delta\) and break ratio \(\beta\), while the other parameter; i.e. the time step \(\text{$\delta $t}\), models the sampling intervals being used during the optimization process. Penalization is a well known technique in optimization that offers the benefit of convexity of the problem. Other beneficial features of penalization is having differential values of cost so that search time for a minimum cost is minimized. In addition, penalization is blocking the parameter from escaping the search field so it is easily localized.

To generate a flexible way to optimize the driving process penalty weights \(p_i\), \(i=1,2,\ldots ,4\) shall be used for the different components of the cost function. The penalty weights are determined according to the importance and major contribution to the cost function. The values are fixed experimentally and allow a broad variation of the optimization process. The cost function uses the method of a least square optimization where different components contribute to the cost function. The following formula comprehends these components:

\begin{equation}
C(\delta ,\beta ,\text{$\delta $t})=\sum _{i=1}^4 p_i\xi _i^2
\label{equation:1}
\end{equation}

\noindent subject to the constraints \(\delta _{\min }\leq \delta \leq \delta _{\max }\), \(\beta _{\min }\leq \beta \leq \beta _{\max }\), and \(\text{$\delta $t}_{\min }\leq \text{$\delta $t}\leq \text{$\delta $t}_{\max }\). The cost function components \(\xi _i\) are measuring the deviation from a predefined target value as:

\begin{equation}
\xi _i=\xi _i^t-\Xi _i\text{      }\text{with}\; \text{    } i=1,2,\ldots ,4\,.
\label{equation:2}
\end{equation}

The values \(\xi _i^t\) are known quantities and shall be fixed in the cost function. As mentioned the \(p_i\) are the penalty weights of the cost function and are also known in the optimization. The objective function is the core of the optimization process and uses the method to minimize the deviation from the ideal line (see above). The ideal line is the reference system to which the dynamical quantities like position, orientation, and velocity are adapted.

The cost function components \(\xi _i\) are assigned to the specific optimization targets. The major target in the optimization process is the alignment of the vehicle to the ideal line. The alignment itself consists of two components \(\xi _1\) and \(\xi _2\) which measure the length of the line element, namely \(\text{$\Delta $s}\), as a quantity of deviation of the actual location of the car to the ideal location using the dynamical position of the vehicle on the road. The second contribution is assigned to the difference of orientation (angular deviation) of the current car dynamics and the ideal orientation of the line using the tangent of the ideal line as a reference. The two quantities depend on the optimization targets the steering angle \(\delta\), the braking ratio \(\beta\), and the time step, \(\text{$\delta $t}\).

As mentioned the first cost component \(\Xi _1\) is assigned to the line element \(\text{$\Delta $s}\). The line element \(\text{$\Delta $s}\) measures the deviation of the position of the car to the ideal position at the current time. Thus, the deviation is determined by the \(x\) and \(y\) coordinates of the vehicle \((x(t),y(t))\) and the ideal line coordinates \(\left(x_l,y_l\right)\). The line element is then:

\begin{equation}
\Xi _1^2=\text{$\Delta $s}^2=\text{$\Delta $x}^2+\text{$\Delta $y}^2=\left(x_l-x(t)\right){}^2+\left(y_l-y(t)\right){}^2
\label{equation:3}
\end{equation}

The optimal solution of the location is reached if \(\text{$\Delta $s}^2=0\) when the vehicle is located on the ideal line. In addition to the location the vehicle has a certain orientation with respect to the ideal line  (see Fig. \cite{geo_setup} for the geometric set up). The direction of the car is determined by the velocity \(\overset{\rightharpoonup }{v}\) with its magnitude \(v=\left. \left. \| \overset{\rightharpoonup }{v}\right. \right. \|\) and its direction as a vector \(\overset{\rightharpoonup }{v}=\left(\dot{x}(t),\dot{y}(t)\right)\). If the current location of the car is used to determine the tangent vector \(\overset{\rightharpoonup }{t}\) to the ideal line there is a complete orientation of the two vectors \(\overset{\rightharpoonup }{t}\) and \(\overset{\rightharpoonup }{v}\) are parallel. This means that according to the relation

\begin{equation}
\overset{\rightharpoonup }{t}.\overset{\rightharpoonup }{v}=\left.\left.\| \overset{\rightharpoonup }{t}\right.\right.\|  \left.\left.\| \overset{\rightharpoonup
}{v}\right.\right.\|  \cos (\theta )
\label{equation:4}
\end{equation}

\noindent we are able to determine the deviation of the orientation as

\begin{equation}
\cos (\theta )=\frac{\overset{\rightharpoonup }{t}.\overset{\rightharpoonup }{v}}{\left.\left.\| \overset{\rightharpoonup }{t}\right.\right.\|  \left.\left.\|
\overset{\rightharpoonup }{v}\right.\right.\| },
\label{eq:5}
\end{equation}

\noindent where \(\theta\) is the alignment angle between the two vectors \(\overset{\rightharpoonup }{v}\) and \(\overset{\rightharpoonup }{t}\). From relation (\ref{eq:5}) it is obvious that the normalized dot product is limited to a finite range \([-1,1]\). The center of this range is known to be related to the orthogonality of the two vectors which is not the target of the optimization. The alignment or anti-alignment of the vectors is given by the end values of the interval. Thus our target is to minimize the quantity 

\begin{equation}
1-\cos (\theta )^2.
\label{equation:6}
\end{equation}

\begin{figure}[t!]
      \includegraphics[scale=0.6]{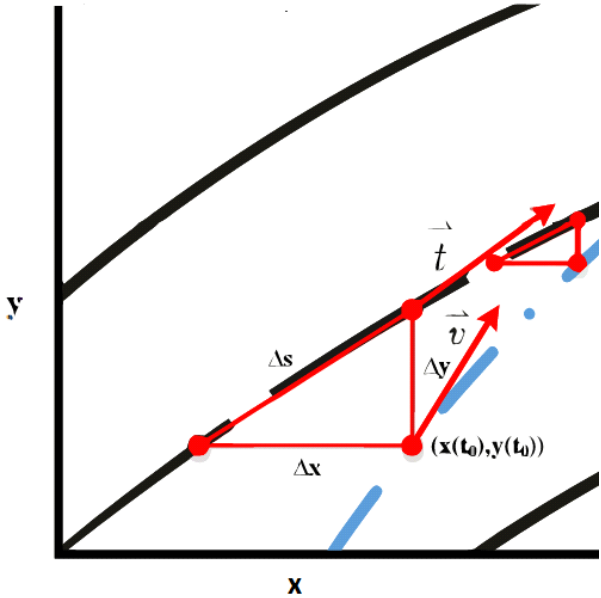}
    \caption{Illustration of how the \(\text{$\Delta $s}\), quantity of deviation, is calculated using \(\text{$\Delta $x}\) and \(\text{$\Delta $y}\) of both the trajectory of the vehicle and center of lane which is the ideal trajectory. Also shown are the tangent vectors of vehicle trajectory and road center lane, \(\overset{\rightharpoonup }{v}\) and \(\overset{\rightharpoonup }{t}\)respectively}
      \label{xa31}
\label{geo_setup}
\end{figure}

The second component \(\Xi _2\) is given by

\begin{equation}
\Xi _2=\cos (\theta )^2=\left(\frac{\overset{\rightharpoonup }{t}.\overset{\rightharpoonup }{v}}{\left.\left.\| \overset{\rightharpoonup }{t}\right.\right.\|
 \left.\left.\| \overset{\rightharpoonup }{v}\right.\right.\| }\right)^2,
\label{equation:7}
\end{equation}

\noindent using the target value \(\xi _2^t=1\) corresponding to an angle \(\theta =0\) or multiples of \(\pi\). The even multiples are related to the alignment while the odd multiples of \(\pi\) are anti-alignments.

Using these two components in an optimization of the vehicle alignment demonstrated that, in principle, a proper position of the vehicle could be found. Moreover, we also observed that above a certain velocity and a certain bending of the road the vehicle went off track. The origin for such a dis-alignment was located in the integration process of the dynamical system. We used in this preliminary optimization approach an equi-distant time sampling to optimize the vehicle alignment. It turned out that such an approach does not take into consideration that a sharper bending of the road will result in an over shooting of the optimized alignment. The solution was to adapt the sampling width \(\text{$\delta $t}\) of the optimization steps to the bending of the road. This optimization of the sampling rate in the optimization process was achieved by limiting the sampling time to a constraint interval 

\begin{equation}
\text{$\delta $t}_{\min }\leq \text{$\delta $t}\leq \text{$\delta $t}_{\max },
\label{equation:8}
\end{equation}

\noindent where the minimal and maximal values were determined by \(\text{$\delta $t}_{\min }=\left.\text{$\delta $t}_1\right/2\) and \(\text{$\delta $t}_{\max }=4 \text{$\delta $t}_1\). The value \(\text{$\delta $t}_1\) belongs to the interval of integration \(\left[t_0,t_E\right]\) where \(t_0\) is the starting point and \(t_E\) the end point in the solution of the Cauchy initial value problem. \(\text{$\delta $t}_1\) which is a mean value of the optimization sampling was found by \(\text{$\delta $t}_1=\left.\left(t_E-t_0\right)\right/N\) where \(N\) is the number of virtual subdivisions of the integration interval. By introducing this third component into the optimization a state where the vehicle was kept on track could be reached. However, if a search for an optimal solution under these conditions was done, it was found that the constrained variation of the braking ration \(\beta\) is not sufficient for the stable operation of the vehicle. It turned out that including the constrained on \(\beta\) as:

\begin{equation}
\beta _{\min }\leq \beta \leq \beta _{\max }
\label{equation:8a}
\end{equation}

\noindent is not stabilizing the dynamic behavior of the car. The minimal and maximal values of \(\beta\) are according to their theoretical definitions given by \(\beta _{\min }=-1\), and \(\beta _{\max }=1\) corresponding to deceleration and acceleration of the vehicle respectively. However, the examination of the dynamical equations previously mentioned demonstrate that the limiting theoretical values are not allowed due to the generation of singularities in the model. The model was adapted in this respect by introducing a limiting parameter \(\epsilon\) which allows us to avoid the singular behavior of the model but keeping the original idea of the breaking ratio. In real computations the values \(\tilde{\beta }_{\min }=\beta _{\min }+\epsilon\) and \(\tilde{\beta }_{\max }=\beta _{\max }-\epsilon\) are used with \(\epsilon\) a small quantity \(\epsilon \approx 10^{-2}\). At this point we note that the used tire model has the disadvantage that the dynamical equations become singular if the limiting values for \(\beta\) are reached.

The observation also was that the velocity in some cases was continuously increased which results to an off road situation. Therefore it was decided in view of the obstacle avoidance conditions to limit the velocity to a given maximal value. The limitation is incorporated by defining

\begin{equation}
\Xi _3=\sqrt{\dot{x}^2+\dot{y}^2}=v
\label{equation:11}
\end{equation}

\noindent and limit the velocity by the target value \(\xi _3^t=v_{\max }\).

Another influence on the vehicle dynamic is generated by the tire model used. The main feature of the tire model is included in the slip angle \(\alpha\) of the model. This quantity is related to the steering angle \(\delta\) for the front wheels and independent of \(\delta\) for rear wheels. To incorporate limitations on the slip angle a penalty function was used in order to avoid an overshoot of the choices of optimization parameters. The standard penalty implementation uses step functions \(u(t)\) to generate the penalty barriers. However, this approach has the disadvantage that the steps introduce discontinuities in the integration process of the Cauchy initial value problem. To avoid discontinuities, and as a result instability of the integration process, continuous penalty functions was introduced by:

\begin{equation}
\label{equation:13}
\begin{split}
\xi _j^2=\lambda -\tanh \kappa  \left(\xi _j^t-\alpha _{j-3}\right)-\tanh \left(\kappa  \left(\xi _j^t+\alpha _{j-3}\right)\right){}^2 \\ \text{
  }\text{with } \text{ } j=4,5,6,7,
\end{split}
\end{equation}

\noindent where \(\lambda\) and \(\kappa\) are constants appropriately chosen. The constraint value on \(\xi _j^t=4{}^{\circ}\) as a fixed angle.

Including all the components discussed above and adding the constraints the optimization problem has the following specific representation:

\begin{align}
C(\delta ,\beta ,\text{$\delta $t})&=\sum _{i=1}^4 p_i\xi _i^2 \nonumber \\
&= \sum _{i=1}^4 p_i \left(\xi _i^t-\Xi _i\right){}^2 
\nonumber \\
&= p_1\left(\left(x_l-x(t)\right){}^2+\left(y_l-y(t)\right){}^2\right){}^2+ \nonumber \\ & p_2\left(1-\left(\frac{\overset{\rightharpoonup
}{t}.\overset{\rightharpoonup }{v}}{\left.\left.\| \overset{\rightharpoonup }{t}\right.\right.\|  \left.\left.\| \overset{\rightharpoonup }{v}\right.\right.\|	
}\right)^2\right){}^2+ 
\nonumber \\
& p_3\left(v_{\max }-\sqrt{\dot{x}^2+\dot{y}^2}\right){}^2+ \nonumber \\
& p_4\sum _{j=4}^7 (\lambda  -\tanh \kappa  \left(\xi _j^t-\alpha _{j-3}\right)) \nonumber \\ & -\tanh
\left(\kappa  \left(\xi _j^t+\alpha _{j-3}\right)\right){}^2 
\end{align}

\noindent subject to \(\delta _{\min }\leq \delta \leq \delta _{\max }\), \(\beta _{\min }\leq \beta \leq \beta _{\max }\), and \(\text{$\delta $t}_{\min }\leq \text{$\delta $t}\leq \text{$\delta $t}_{\max }\). Thus the final problem is to find the minimum of the cost function

\hfill \break

\begin{equation}
\underset{\delta ,\beta ,\text{$\delta $t}}{\text{Min}}(C(\delta ,\beta ,\text{$\delta $t}))
\label{equation:15}
\end{equation}

\begin{multline}
\noindent \text{under} \text{ the} \text{ constraint }  \delta _{\min }\leq \delta \leq \delta _{\max }, \beta _{\min }\leq \beta \leq \beta _{\max },\\  \text{ and }\text{$\delta $t}_{\min }\leq \text{$\delta $t}\leq \text{$\delta $t}_{\max }.
\label{equation:16}
\end{multline}



\begin{figure}[t!]
     \centering
      \includegraphics[scale = 0.69]{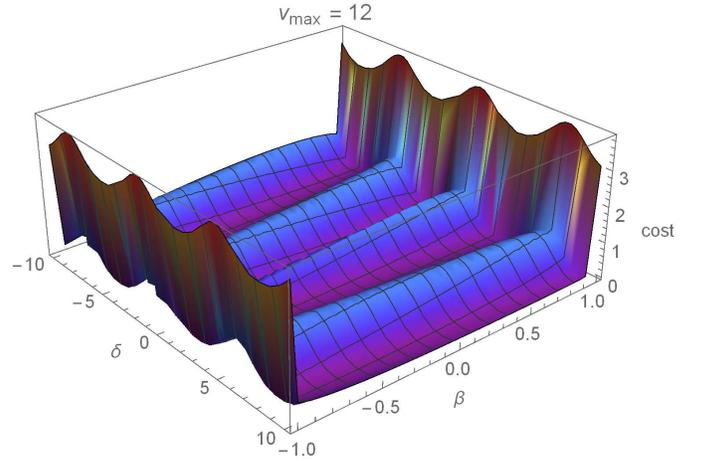}
    \caption{Sample for the cost function (\(C\) = cost) for a straight inclined lane as road model within time interval 0.05s at velocity of 12m/s using the feature of penalization by adding penalty weights of 0.5 to each component of the cost function.}
      \label{Figure-2}

\end{figure}

\subsubsection{Optimization Parameters}

The field-of-search parameters are determined as two types. One type is concerned about optimizing the correspondent driver's steering angle and vehicle's braking ratio to a minimum objective function cost while the other is correlated to optimizing the time samples. Listed below are two different types of optimization parameter that were used throughout this work.

\begin{figure}[t!] 
  \includegraphics[height=0.2\textheight,width=0.45\textwidth]{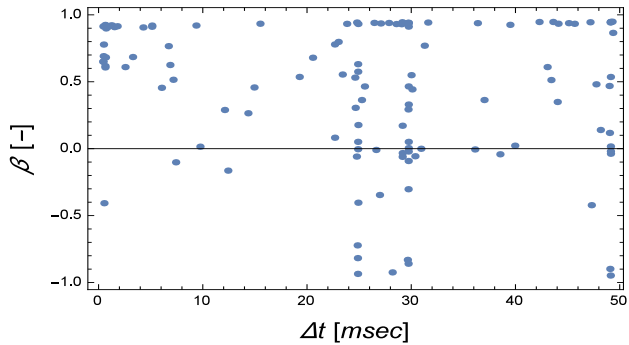}
\caption{An example of time samples used during optimization done in a scenario of a rough obstacle road model simulation. The scenario is shown in a test illustrated in Fig. \ref{qqqex}. Shown is the braking ratio optimized values and the different time samples which were used throughout the control process.}
\label{qqq}
\end{figure}

\begin{enumerate}

\item[A.] \textit{Vehicle's Parameters}

The optimization of the parameters were done at each time-sample by finding the most optimum-minimum cost then relating the $\alpha$ and $\beta$ to it at $t_{curr}$.

\item[B.] \textit{Time Sampling Parameter}

An internal optimization is done during the optimization process in order to search for the current time sample size that increases the possibility to find a minimum cost. Sampling time intervals into equidistant intervals was found to not have an optimum solution according to \cite{c}, \cite{d}, in contrast to optimizing time samples and searching for different time sample at every iteration that was found to be more-likely to find an optimum solution and converge eventually. Shown in Fig. \ref{qqq} is an example of an optimization using time samples throughout the whole process. The figure also shows that different time intervals where used in order to avoid a possible overshoot in the control process of the vehicle.
\end{enumerate}

\begin{figure}[t!]
   \centering
    \includegraphics[scale = 0.5]{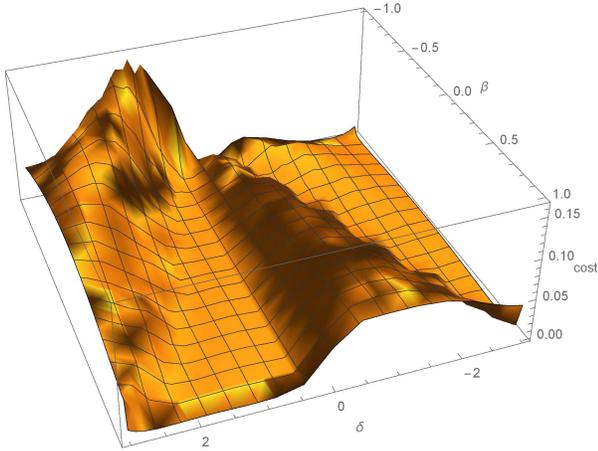}
  \caption{Sample for the cost function (\(C\) = cost) for a straight inclined lane as road model with no penalty weights implemented at velocity of 12m/s.}
      \label{Figure-3}
\end{figure}

\begin{figure}[t!] 
  \includegraphics[height=0.2\textheight,width=0.45\textwidth]{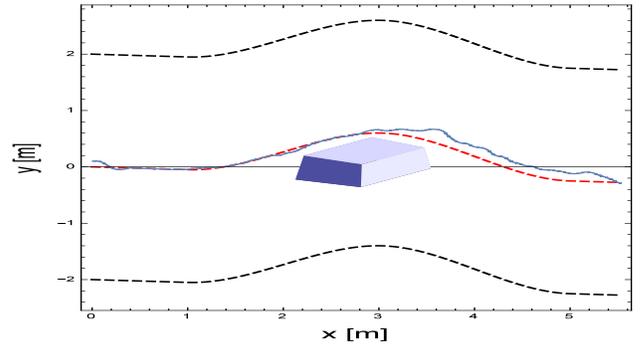}
\caption{A tested scenario of a vehicle approaching a rough obstacle, $\epsilon = 0.5$, on the road where the optimization and control algorithm succeeded to follow the ideal road and avoid the presented rough obstacle simultaneously.}
\label{qqqex}
\end{figure}

\subsubsection{Constraints and Penalization}

The developed objective function had to follow a constraint definition and penalty for constraints violation in order to insure vehicle safety. Slip angle constraint is the insurance that the vehicle never slips out of the road which leads to rotation uncontrollably of the vehicle and entering an unsafe state. Such constraint is defined to be in the range of \([-4,4]\) deg according to \cite{20}. Another constraint that is added is that $\beta$ can not exceed the boundaries of \([-1,1]\). This constraint is added for mathematical safety more than a vehicle physical safety, i.e. to insure elimination of infinite results as in dividing by zero. 
 
Penalty function is a value that is multiplied into the cost function in order to be able to localize these parameters and find a conversion field where the parameters converge to a certain value. In Fig. \ref{Figure-2}  there was no penalty implemented and that is why the optimizer could not find a solution or converging values for the parameters. In Fig. \ref{Figure-3} , penalty function was implemented in order to lock the parameters into barriers so the parameters could converge. Penalty weight for constraints violation is determined by trial and error of random values according to the system behavior. The shown cost function was plotted for an interval of $0.05s$ at a velocity of $12m/s$.



\textbf{Proposed Optimization Algorithm}

\begin{algorithm}
\label{fig:optim}
\begin{algorithmic}[1]
\Procedure INPUTS:$[t_{start}, t_{end}, init, R_{line}, V_{max} 
, P_{wg}]$
\State $t_s = t_{start}$
\State $t_e = t_{end}$
\State $\delta t_{max} = \frac{t_e-t_s}{div}$
\State $\delta t = \delta t_{max}$
\State $t_e = t_s + \delta t$
\While $(t_e < t_{end})$
\State FindMinimum[\(C\), $init$, $t_s$, $t_e$] for $\delta$ and $\beta_r$ and $\delta t$
\State Subject to constraints $\delta_{min} < \delta_{new} < \delta_{max}$, $\beta_{min} < \beta_{new} < \beta_{max}$, $0.05 	\delta t_{max} < \delta t_{new} < \delta t_{max}$
\State $\delta = \delta_{new}$
\State $\beta_r = \beta_{new}$
\State Append $\delta$ to vector of solutions
\State Append $\beta_r$ to vector of solutions
\State $t_s = t_e$
\State $t_e = t_s + \delta t_{new}$
\EndWhile
\State \textbf{return} $\delta$ and $\beta_r$ vectors \Comment{Interpolation used}
\EndProcedure
\end{algorithmic}
\caption{Optimization algorithm for collecting the data of each iteration and interpolate it into a vector that can be visualized after the process is done.}
\end{algorithm}

The algorithm is illustrated above, where $t_{start}$ and $t_{end}$ are the limits of the whole simulation interval, $t_s$ and $t_e$ are the limits of each optimization iteration based on the discrete sampling. $\delta_{new}$ and $\beta_{new}$ are the new optimized steering angle and braking ratio for the current iteration that corresponds to a minimum \(C\) which is the cost function. $init$ are the initial values for $x$, $\dot{x}$, $y$, $\dot{y}$, $\psi$, and $\dot{\psi}$. $div$ denotes the number of divisions that the whole interval should be sampled to. $R_{line}$ denotes the ideal line that the optimization should refer to when calculating the minimum cost function. $V_{max}$ and $P_{wg}$ denotes the maximum allowed velocity and constraints violation penalty weights respectively.

\section{Results}
The optimization algorithm presented in the previous section is tested and validated along vehicle model, and road model previously discussed. The driver model is overridden by the optimization and control throughout the simulation. Test cases include low speed for straight lane model. Avoiding an obstacle is also tested and validated for low speed vehicle in this section. High speed scenarios are tested for inclined lane model.


The vehicle parameters chosen for the test and validation are shown in Table \ref{table:1}. The parameters are used in order to simulate the reaction of a real vehicle in the simulation carried out through this work. These parameter represent the dimensions and forces of a commercial vehicle.

\begin{table}[H]
\def\arraystretch{1}
\caption{Vehicle parameters chosen for the simulation test.}
\centering
\begin{tabular}{p{2cm}|p{2cm}|p{2cm}}

 \hline
 Parameter & Value & Units\\
 \hline
 $w_t$ & $1.63$ & $m$\\
 $l_f$ & $1.43$    &  $m$\\
 $l_r$ & $1.47$ & $m$\\
 $\mu$ & $1$  & -\\
 $C_\alpha$&   $80,000$  & -\\
 $m$ & $2050$  & $kg$  \\
 $Fz_{1,2}$& $26,719$  & $N$\\
 $Fz_{3,4}$& $21,295$  & $N$\\
 $J_z$& $3344$  & $kg.m^2$\\
 \hline

\end{tabular}

\label{table:1}
\end{table}

The optimization and control parameters chosen throughout the simulation will be mentioned for every test case in the following sections. The optimization parameters include simulation time, penalty weights, initial values, road length, and obstacle width. Some parameters were globalized for all test cases such as road length, steering angle limits, discrete intervals, braking ratio limits, and slip angle limits. The global parameters are shown in Table \ref{table:2}. Listed below are the different test scenarios which were simulated and verified.

\begin{table}[H]
\def\arraystretch{1}
\caption{Optimization and design parameters chosen for the simulation test.}
\centering
\begin{tabular}{p{2cm}|p{2cm}|p{2cm}}

 \hline
 Parameter & Value & Units\\
 \hline
 $\text{Road Length}$ & $100$ & $m$\\
 $\delta _{\min }$ &  $- 0.06$    &  $rad$\\
 $\delta _{\max }$ & $0.06$ & $rad$\\
 $\delta t_{\max }$ &  $\frac{t_{end} - t_{start}}{\text{No. of divisions}}$  & $sec$\\
 $\delta t_{\min }$&   $0.05 * \delta t_{\max }$   & $sec$\\
 $\beta_{min}$ & $-1$  & -  \\
 $\beta_{max}$& $1$  & -\\
 $\alpha_{min}$& $-4/2\pi$  & $rad$\\
 $\alpha_{max}$& $4/2\pi$  & $rad$\\
 \hline

\end{tabular}

\label{table:2} 
\end{table}

\subsubsection{Maintaining a parabolic road at a velocity 18 m/s}

A tested scenario for the vehicle at an initial velocity of 18 m/s was done using a road model of a parabolic surface. The goal of the control algorithm is to maintain the ideal road and keeping cost function at minimal value. The optimization and design parameters chosen for this test  is shown in Table \ref{table:3}.

\begin{table}[H]
\def\arraystretch{1.3}
\caption{Optimization parameters chosen for a parabolic road geometry at a velocity of 18 m/s.}

\centering
\begin{tabular}{p{1cm}|p{1cm}|p{1cm}||p{1cm}|p{1cm}|p{1cm}}

 \hline
 Parameter & Value & Units & Parameter & Value & Units\\
 \hline
 $x$ & $0.01 $ & $m$ & $\dot{x}$ & $18 $& $\frac{m}{s}$  \\
 $y$ & $0.1 $ & $m$ & $\dot{y}$ & $0.5 $& $\frac{m}{s}$  \\
 $N_d$ & $250$ & $-$ & $v_{max}$ & $2.5 $& $\frac{m}{s}$  \\
 $p_1$ & $0.45$ & $-$ & $p_2$ & $0.45$& $-$  \\
 $p_3$ & $0.1$ & $-$ & $p_4$ & $1$& $-$  \\
 $\epsilon$ & $0.7 $  & $m$ & $t_{sim}$ & $2 $& $sec$  \\
 \hline
\end{tabular}
\label{table:3}
\end{table}

\begin{figure}[t!] 
\begin{subfigure}{.45\textwidth}
  \includegraphics[height=0.2\textheight,width=\textwidth]{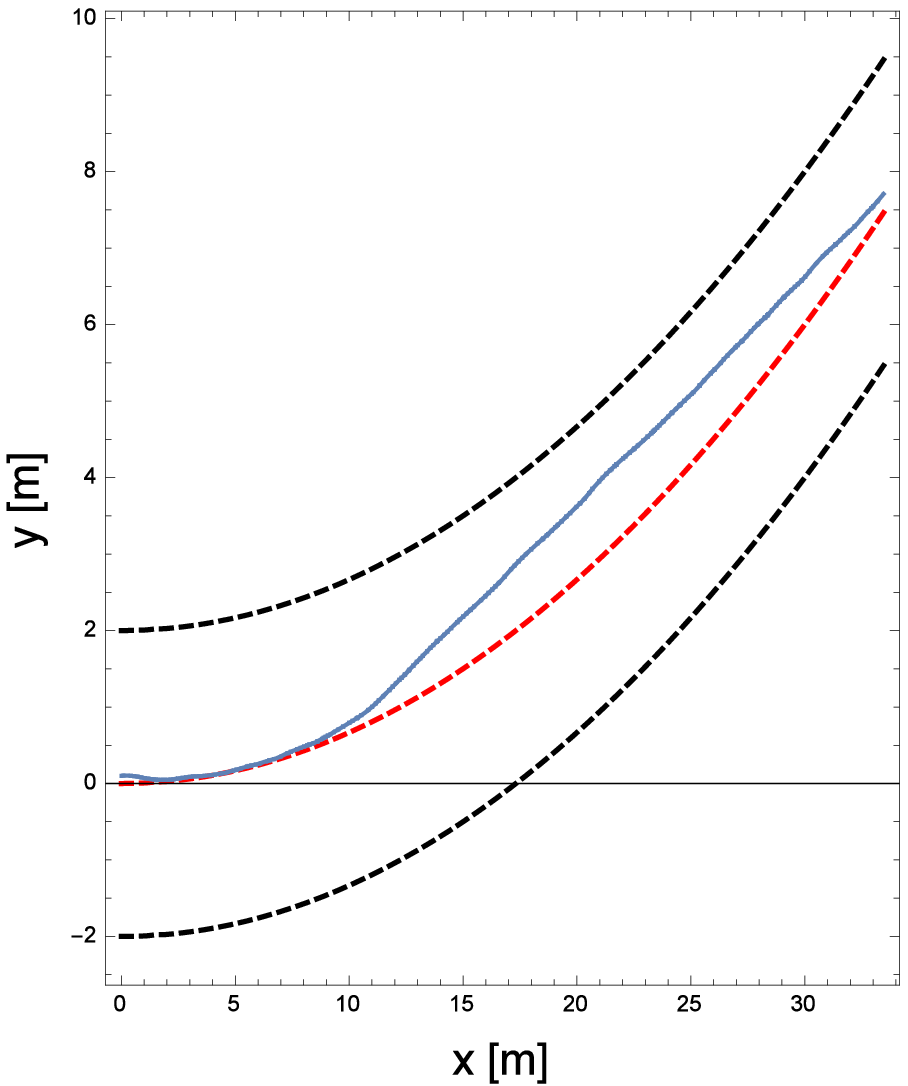}
  \caption{ }
  \label{}
\end{subfigure}
\hspace{\fill}
\begin{subfigure}{.45\textwidth}
  \includegraphics[height=0.2\textheight,width=\textwidth]{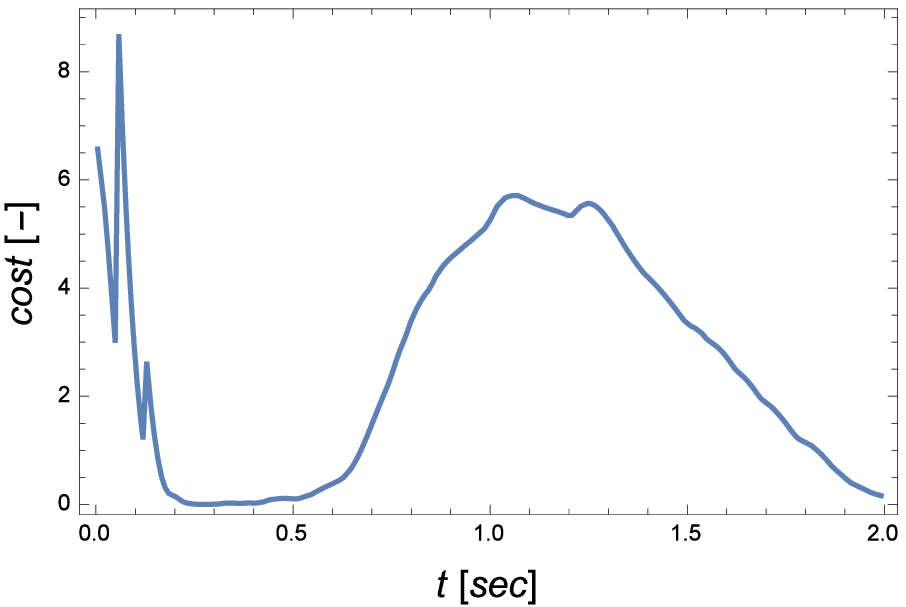}
  \caption{ }
  \label{}
\end{subfigure}
\caption{Simulation results. In these plots, we show the behavior of the control algorithm on maintaining the vehicle to follow the ideal road geometry. The plots illustrates the road geometry which represent an input to the control algorithm, and the vehicle’s performance which represent the output of the control algorithm. (a) The parametric plot result of vehicle positions done by the optimization and control algorithm at an initial input velocity of $18 m/s$. The centered red dashed plot represents the ideal road while the sided black dashed plots represent left and right lane margin. (b) The plot indicates the change of the cost function during the process. }
\label{q1}
\end{figure}

\begin{figure}[t!] 
\begin{subfigure}{.45\textwidth}
  \includegraphics[height=0.15\textheight,width=\textwidth]{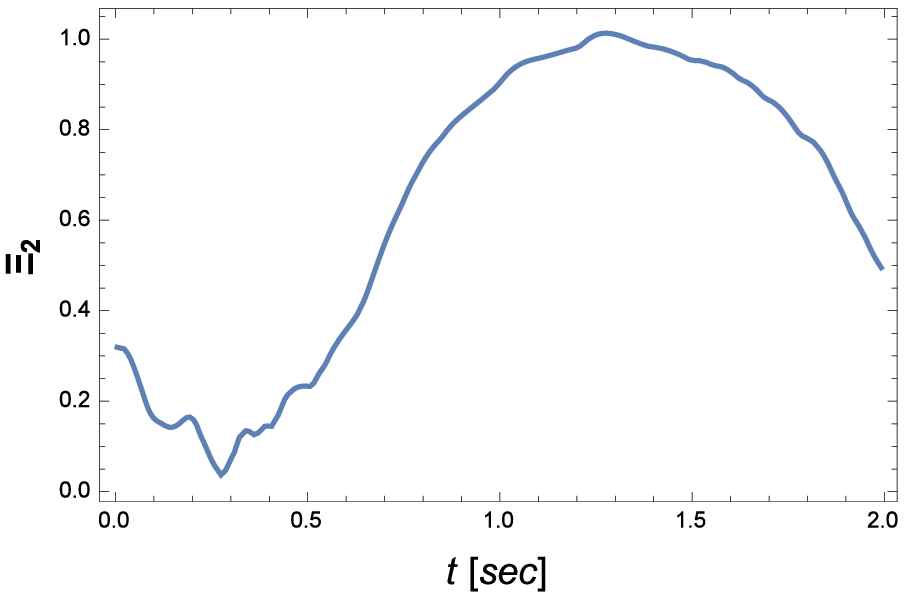}
  \caption{ }
  \label{fig:sub2}
\end{subfigure}
\hspace{\fill}
\begin{subfigure}{.45\textwidth}
  \includegraphics[height=0.15\textheight,width=\textwidth]{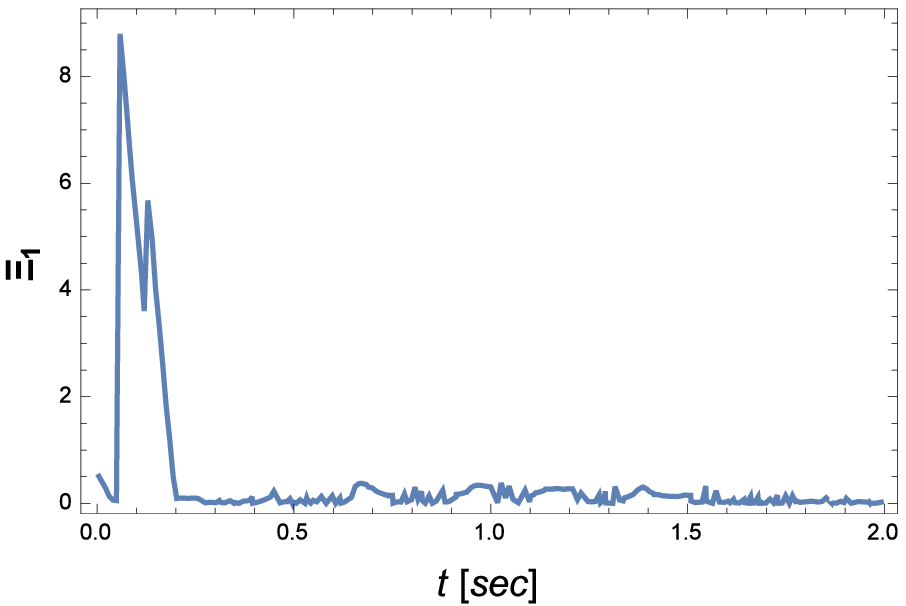}
  \caption{ }
  \label{fig:subs}
\end{subfigure}
\caption{ Simulation results. In these plots, the mathematical approach is analyzed by observing the change of two of the cost function components, orientation error and deviation quantity for the tested scenario at a velocity of 18 m/s.  (a) Orientation error of the vehicle with respect to the presented road geometry is indicated in this plot. (b) The vehicle's deviation quantity plot is shown for the illustration of the vehicle's performance}
\label{q2}
\end{figure}

\begin{figure}[t!] 
\begin{subfigure}{.45\textwidth}
  \includegraphics[height=0.15\textheight,width=\textwidth]{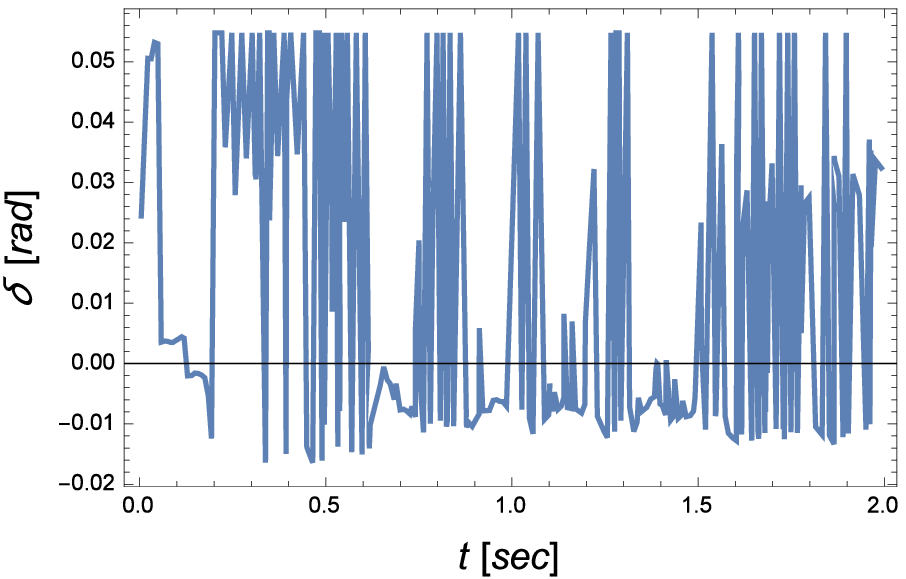}
  \caption{ }
  \label{fig:sub1}
\end{subfigure}
\hspace{\fill}
\begin{subfigure}{.45\textwidth}
  \includegraphics[height=0.15\textheight,width=\textwidth]{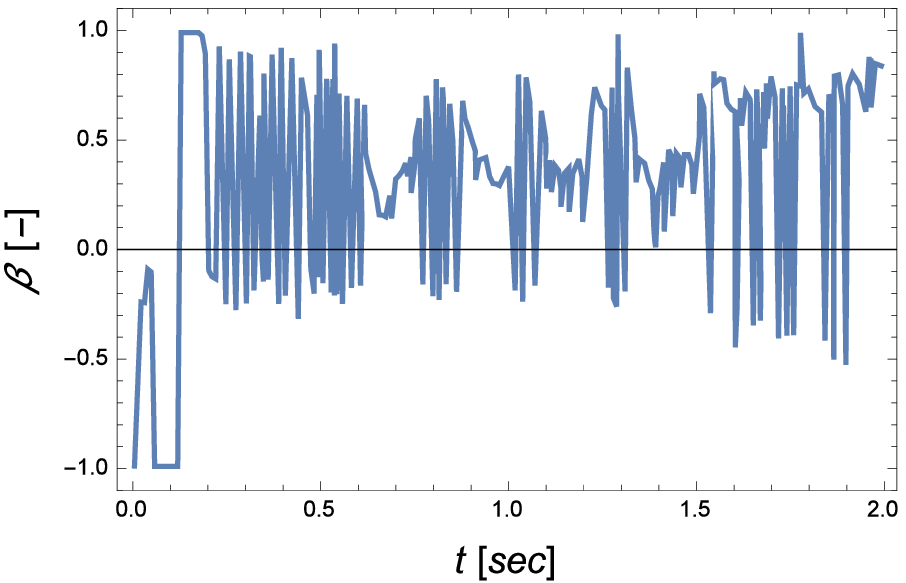}
  \caption{ }
  \label{fig:sub2}
\end{subfigure}
\caption{ Simulation results. These plots capture the steering angle and braking ratio readings during the optimization process of a vehicle following a steep parabolic road geometry at a velocity of 18 m/s. These simulation results are compatible with a simulator. (a) The control algorithm is manipulating the steering angle in order to follow the input road geometry. (b) Braking ratio value is changing rapidly in order to maintain the vehicle within the maximal velocity limits.}
\label{q3}
\end{figure}

\begin{figure}[t!] 
\begin{subfigure}{.45\textwidth}
  \includegraphics[height=0.15\textheight,width=\textwidth]{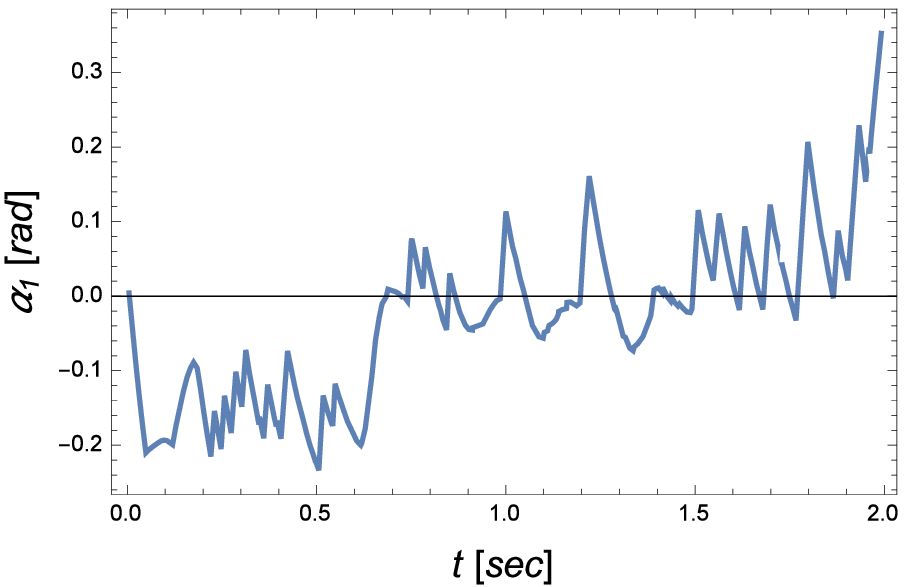}
  \caption{ }
  \label{fig:sub1}
\end{subfigure}
\hspace{\fill}
\begin{subfigure}{.45\textwidth}
  \includegraphics[height=0.15\textheight,width=\textwidth]{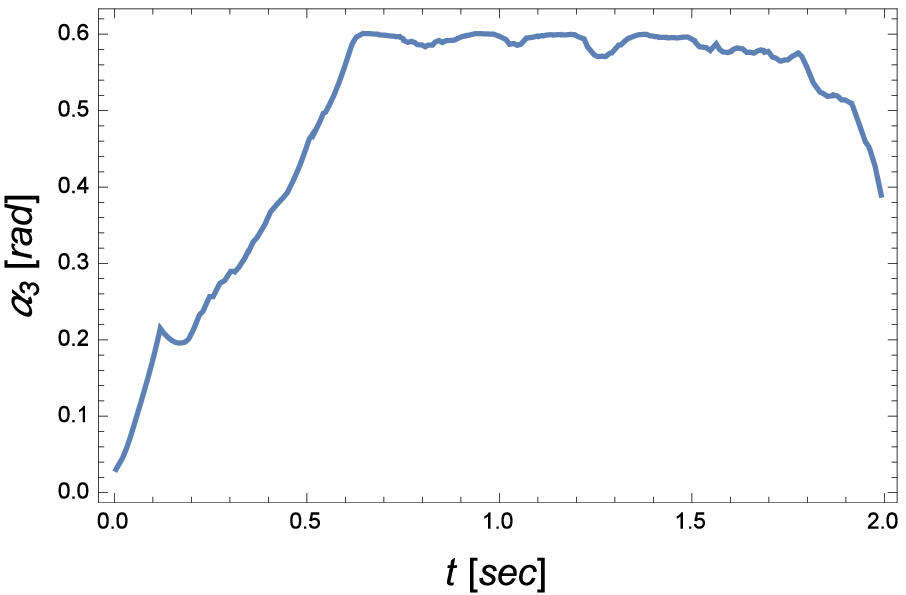}
  \caption{ }
  \label{fig:sub2}
\end{subfigure}
\caption{Plot of tires states showing the change in slip angles during the optimization and control process of a vehicle following a parabolic road geometry at velocity of 18 m/s. (a) Plot of front tires slip angle. (b) Plot of rear tires slip angle.}
\label{q4}
\end{figure}

It is noted in Fig. \ref{q1}(a), which shows the vehicle's performance, that the control algorithm could successfully maintain the vehicle close to the ideal road geometry with a minimal error that is shown and analyzed in Fig. \ref{q2}, \ref{q3}, and \ref{q4}. Fig. \ref{q1}(b) shows the cost function plot results throughout the duration of the control and optimization process. The overall cost function value throughout the simulation indicates that the optimization was successful at keeping the vehicle positions within the safe limits of the ideal road to-follow.

Shown in Fig. \ref{q2}(a) is the deviation from the optimum position, the center of lane, of a road model with a presented parabolic road surface throughout the simulation time. Also in Fig. \ref{q2}(b) the orientation error for the vehicle during the control process is shown.

In Fig. \ref{q3}(a), we note that the steering angle was maintained within the constraint value which was preset during the optimization process. Fig. \ref{q3}(b) shows the interpolated data of the braking ratio throughout the whole simulation process. The data shows that the braking ratio was applied frequently and that was to maintain the maximum velocity that was specified in the design parameters. The braking ratio is changing rapidly to maintain the velocity specified at a steep parabolic road which was challenging. We note from these results that the constraints set for steering angle, braking ratio, and maximal velocity were successfully met throughout this test.

Fig. \ref{q4}(a) shows the slip angle data of the front tires of the vehicle on a road model with a presented rough obstacle. The figure shows that the slip angle limits were not exceeded. The optimization algorithm has succeeded to meet the slip angle constraint. The optimization control shows a good performance in controlling the vehicle to follow the optimum path. Fig. \ref{q4}(b) shows the slip angle data for the rear tires at each time instant. It is shown that the slip angle is exactly between the maximum and minimum limits. The optimization algorithm has succeeded to meet the set constraint for the rear slip angles.

The simulation result and data previously discussed shows that the proposed optimization and control algorithm could successfully control the vehicle to maintain the presented parabolic road geometry without violating any of the presented constraints for the vehicle motion. Moreover, the algorithm was successful at meeting the steering angle, and braking ratio constraint as shown in Fig. \ref{q3}. The optimization control was also able to successfully maintain the vehicle on the ideal road while minimizing the cost function to an optimum value as shown in Fig. \ref{q2}(b) and \ref{q1}(a). The slip angle constraints were met successfully for front and rear tires as shown in Fig. \ref{q4}. A summary of the evaluation done for the scenario test is shown in Table \ref{s6}.

\begin{table}[h!]
\def\arraystretch{1.5}
\centering
\begin{tabular}{||p{2cm}|p{2cm}|p{2cm}||}

 \hline
 Test & Value &Evaluation\\
 \hline
 $\text{Vehicle Motion}$ & - & $\text{Success}$\\
 $\text{Steering Angle}$ & $\in [-0.018, 0.052]$ & $\text{Success}$\\
 $\text{Braking Ratio}$ & $\in [-1,1]$ & $\text{Success}$\\
 $\text{Front Slip Angle}$ & $\in [-0.2, 0.3]$ & $\text{Success}$\\
 $\text{Rear Slip Angle}$& $\in [0, 0.6]$ & $\text{Success}$ \\
 \hline

\end{tabular}

\caption{Scenario evaluation for a parabolic road geometry at a velocity of 18 m/s.}
\label{s6}
\end{table}

\subsubsection{Approaching an obstacle at velocity of 10 m/s}
Another scenario was tested such that a parabolic road was presented in addition to an obstacle. The optimization and control algorithm is tested by observing the vehicle's performance towards following the ideal path and avoiding the presented obstacle. The optimization and design parameters chosen for this test case is shown in Table \ref{table:4}.

\begin{table}[H]
\def\arraystretch{1.3}
\caption{Optimization parameters chosen for a vehicle velocity of 10 m/s approaching an obstacle.}
\centering
\begin{tabular}{p{1cm}|p{1cm}|p{1cm}||p{1cm}|p{1cm}|p{1cm}}

 \hline
 Parameter & Value & Units & Parameter & Value & Units\\
 \hline
 $x$ & $0.01 $ & $m$ & $\dot{x}$ & $10 $& $\frac{m}{s}$  \\
 $y$ & $-1 $ & $m$ & $\dot{y}$ & $1 $& $\frac{m}{s}$  \\
 $N_d$ & $50$ & $-$ & $v_{max}$ & $9.5 $& $\frac{m}{s}$  \\
 $p_1$ & $0.75$ & $-$ & $p_2$ & $0.45$& $-$  \\
 $p_3$ & $0.1$ & $-$ & $p_4$ & $0.1$& $-$  \\
 $\epsilon$ & $0.1$  & $m$ & $t_{sim}$ & $2 $& $sec$  \\
 \hline
\end{tabular}
\label{table:4}
\end{table}

\begin{figure}[t!] 
\begin{subfigure}{.45\textwidth}
  \includegraphics[height=0.25\textheight,width=\textwidth]{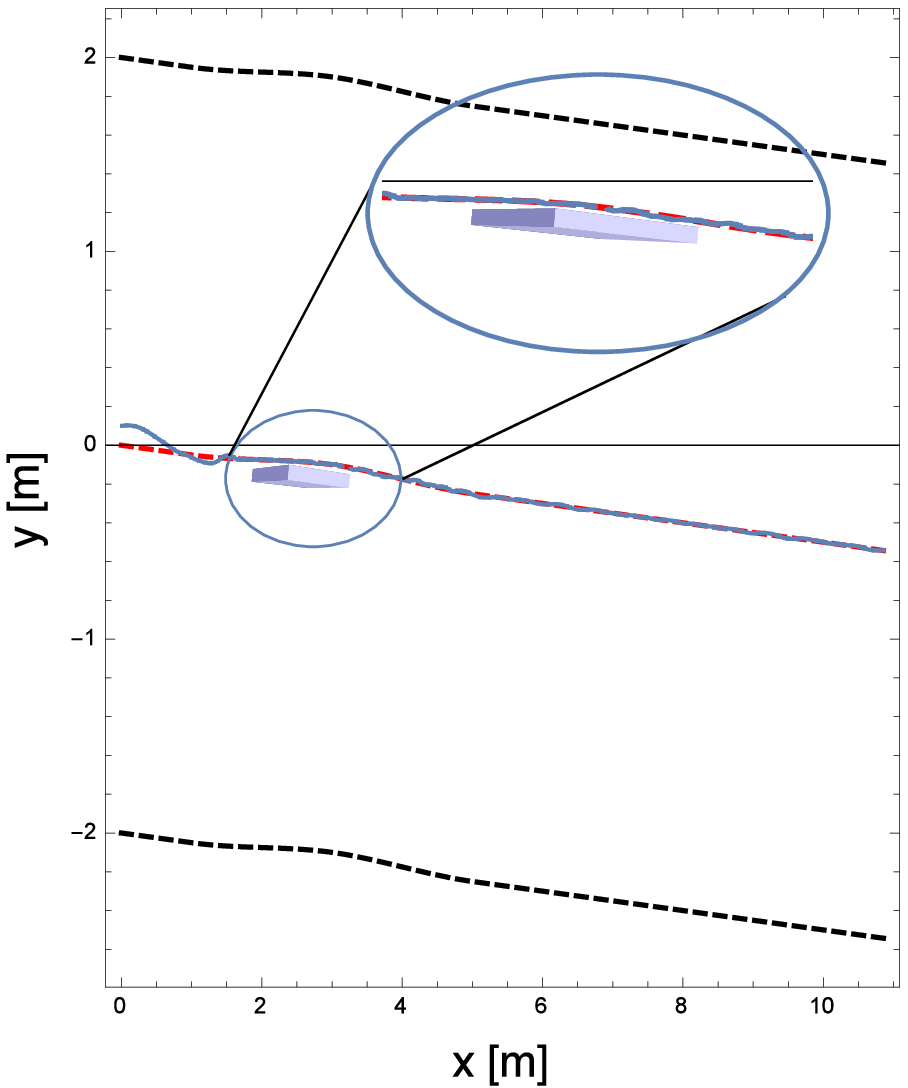}
  \caption{ }
  \label{fig:sub1}
\end{subfigure}
\hspace{\fill}
\begin{subfigure}{.45\textwidth}
  \includegraphics[height=0.25\textheight,width=\textwidth]{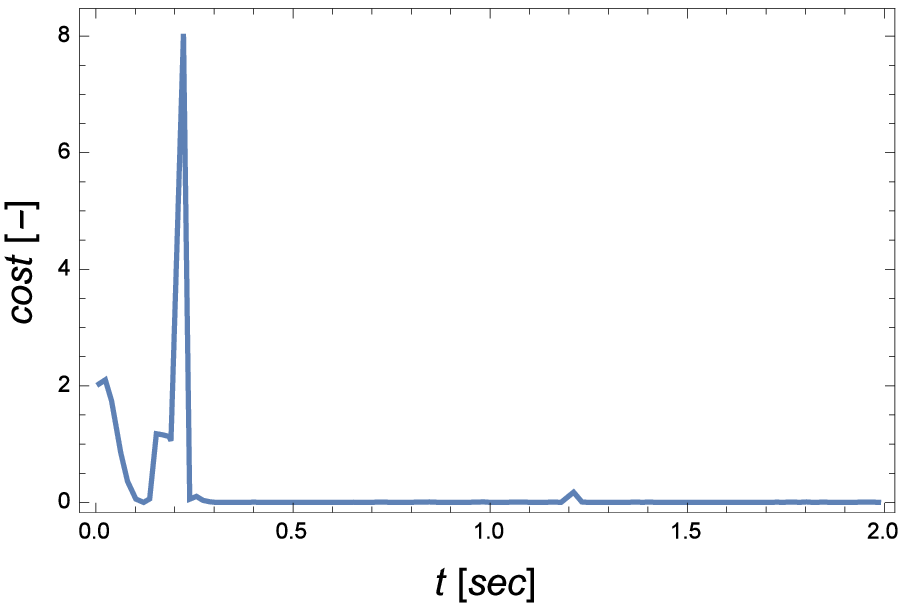}
  \caption{ }
  \label{fig:sub2}
\end{subfigure}
\caption{Simulation results. In these plots, we show the behavior of the control algorithm on maintaining the vehicle to follow the ideal road geometry presented with a smooth obstacle of width $\epsilon = 0.1$. The plots illustrates the road geometry which represent an input to the control algorithm, and the vehicles performance which represent the output of the control algorithm. (a) The parametric plot result of vehicle positions done by the optimization and control algorithm at an initial input velocity of $10m/s$. The centered red dashed plot represents the ideal road while the sided black dashed plots represent left and right lane margin. (b) The plot indicates the change of the cost function during the process. Vehicle’s initial conditions were chosen in a way that a miss-alignment between vehicle’s motion and road’s geometry occurred which led to a transient amplitude. It is observed that this transient behavior is vanished throughout the rest of the simulation after the controller has adjusted vehicle’s position to road’s geometry.}
\label{e1}
\end{figure}

Shown in Fig. \ref{e1} is the vehicle performance on the road model presented and the relative cost function plot during the control process. The road geometry presented for this scenario's evaluation test, as discussed in Section II. B., has a parabolic geometry with a presented smooth obstacle which is required to be avoided by the optimization and control algorithm. It is noted in Fig. \ref{e1}(a) that the controller could successfully lead the vehicle to avoid the obstacle and keeping the track of the ideal road. 

\begin{figure}[t!] 
\begin{subfigure}{.45\textwidth}
  \includegraphics[height=0.15\textheight,width=\textwidth]{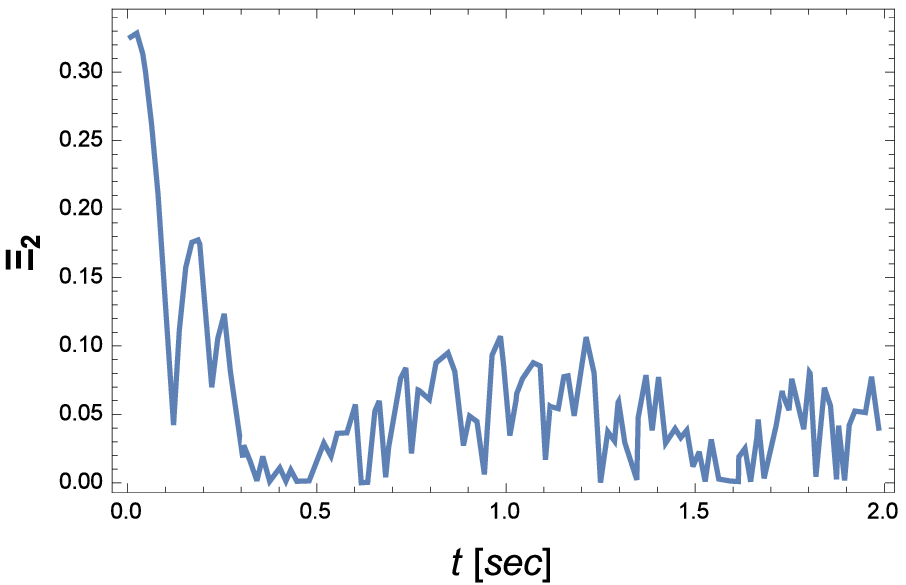}
  \caption{ }
  \label{fig:sub2}
\end{subfigure}
\hspace{\fill}
\begin{subfigure}{.45\textwidth}
  \includegraphics[height=0.15\textheight,width=\textwidth]{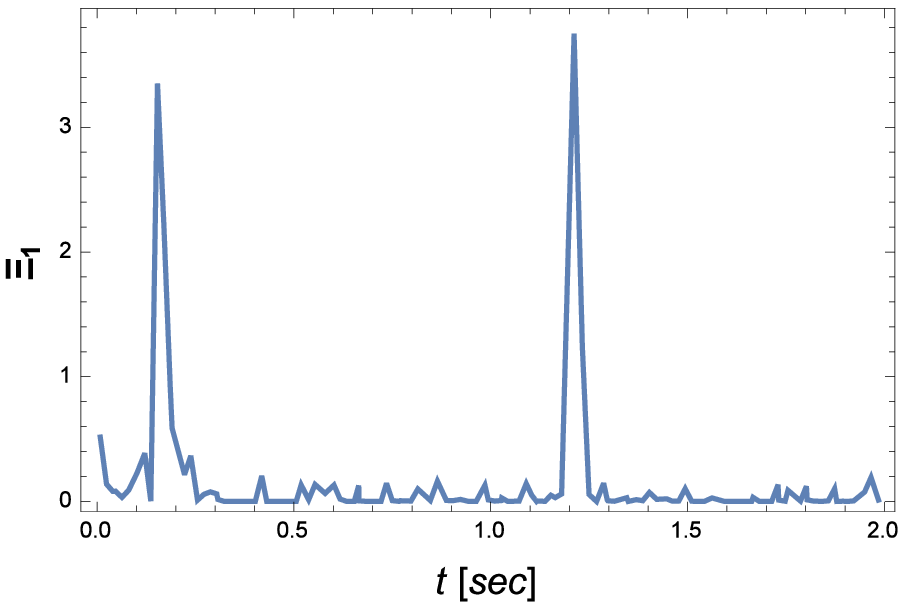}
  \caption{ }
  \label{fig:subs}
\end{subfigure}
\caption{ Simulation results. In these plots, the mathematical approach is analyzed for the tested scenario of a vehicle controlled to follow the ideal path and avoid an obstacle at a velocity of 10 m/s. (a) Orientation error of the vehicle with respect to the presented road geometry is indicated in this plot. (b) The vehicle's deviation quantity plot is shown for the illustration of the vehicle's performance}
\label{e2}
\end{figure}

We note in Fig. \ref{e1}(b) that the cost function converges to zero-value with a minimal error that is analyzed in Fig. \ref{e2}(a) and Fig. \ref{e2}(b) from a mathematical point of view. Fig. \ref{e2}(a) shows the deviation from the optimum path which is presented as a parabolic road geometry including an obstacle. Moreover, Fig. \ref{e2}(b) shows the orientation error for the vehicle during the simulation process. It is important to view these plots in order to understand the mathematical behavior of the optimization and control algorithm while controlling the vehicle to follow the road and avoid the presented obstacle. We note that the orientation error and deviation quantity plots show that the vehicle was in parallel and oriented with the ideal road throughout the simulation process.

\begin{figure}[t!] 
\begin{subfigure}{.45\textwidth}
  \includegraphics[height=0.15\textheight,width=\textwidth]{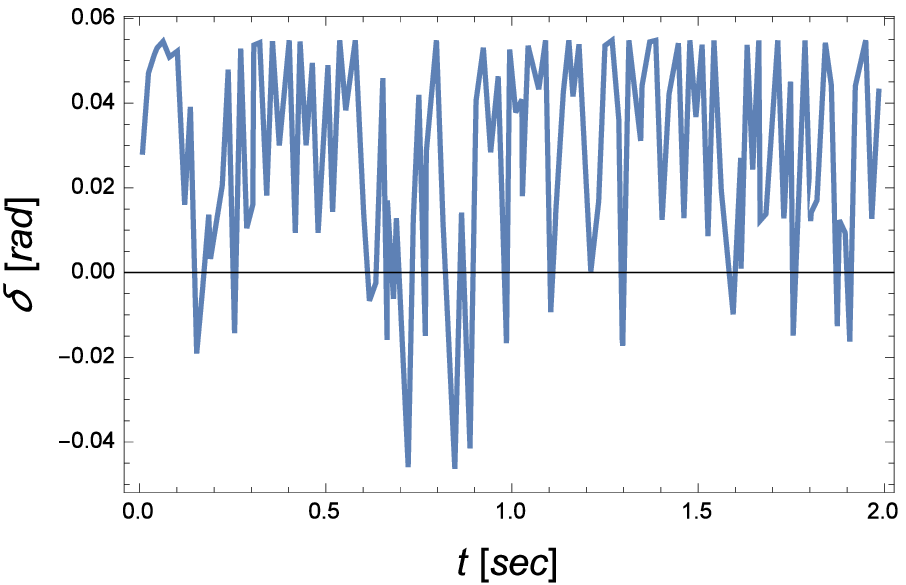}
  \caption{ }
  \label{fig:sub1}
\end{subfigure}
\hspace{\fill}
\begin{subfigure}{.45\textwidth}
  \includegraphics[height=0.15\textheight,width=\textwidth]{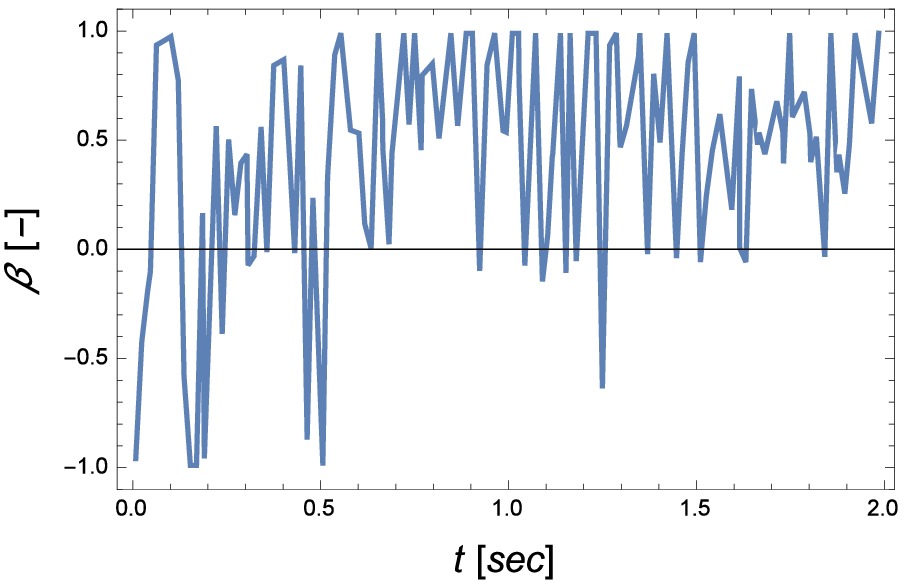}
  \caption{ }
  \label{fig:sub2}
\end{subfigure}
\caption{These plots capture the steering angle and braking ratio readings during the optimization process of a vehicle following a parabolic road geometry and avoiding an obstacle at a velocity of 10 m/s. These simulation results are compatible with a simulator. (a) The control algorithm is manipulating the steering angle in order to follow the input road geometry. (b) Braking ratio value is changing rapidly in order to maintain the vehicle within the maximal velocity limits.}
\label{e3}
\end{figure}

We note in Fig. \ref{e3}(a) that the steering angle data that mimics the driver behavior didn't exceed the constraint which was set in the control algorithm. Also we noted in Fig. \ref{e3}(b) that the interpolated data of the braking ratio throughout the whole simulation process shows that the minimum and maximum limits of $\beta$ was not exceeded and thus the optimization constraint for this part was not violated.

\begin{figure}[t!] 
\begin{subfigure}{.45\textwidth}
  \includegraphics[height=0.15\textheight,width=\textwidth]{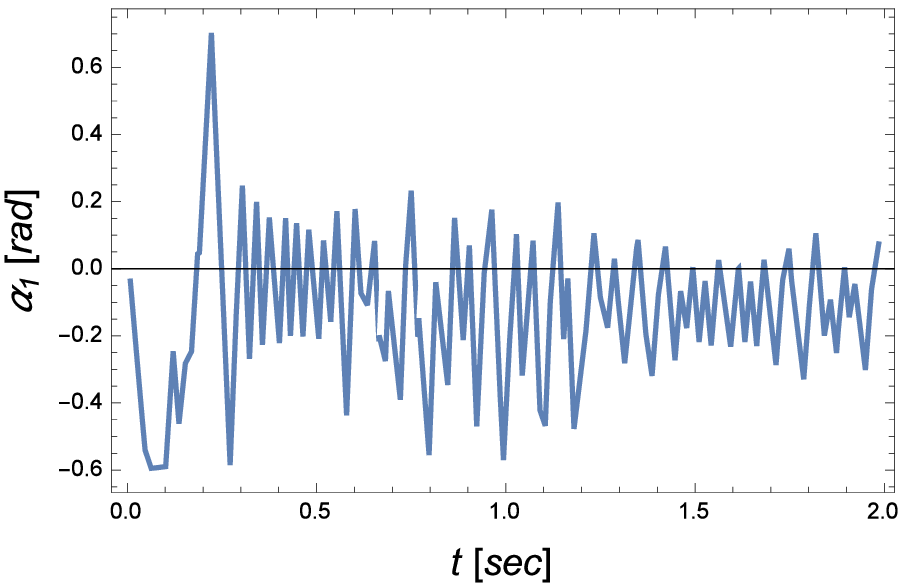}
  \caption{ }
  \label{fig:sub1}
\end{subfigure}
\hspace{\fill}
\begin{subfigure}{.45\textwidth}
  \includegraphics[height=0.15\textheight,width=\textwidth]{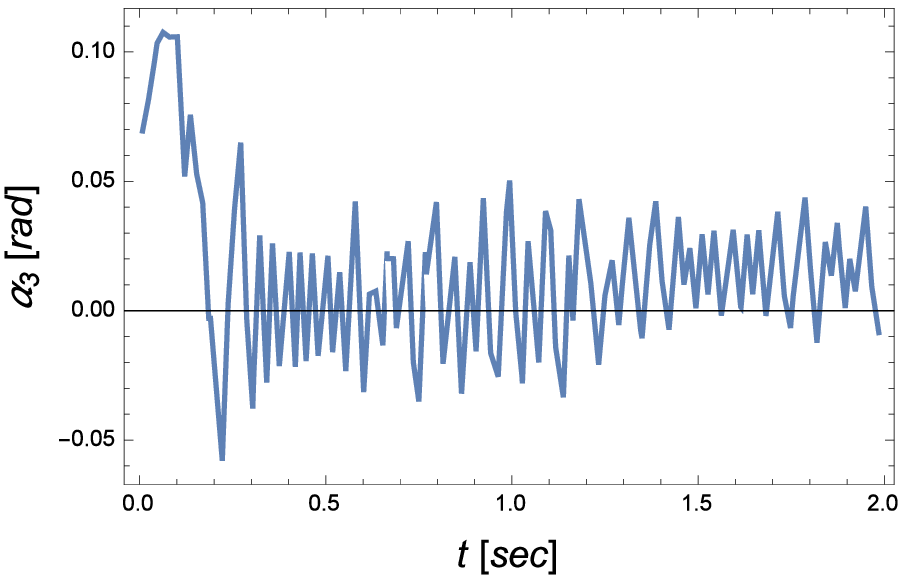}
  \caption{ }
  \label{fig:sub2}
\end{subfigure}
\caption{Plot of tires states showing the change in slip angles during the optimization and control process of a vehicle following a parabolic road geometry and avoiding an obstacle at velocity of 10 m/s. (a) Plot of front tires slip angle. (b) Plot of rear tires slip angle.}
\label{e4}
\end{figure}

Shown in Fig. \ref{e4}(a) is the slip angle data of the front tires of the vehicle on parabolic road model with a presented obstacle. The figure shows that the slip angle limits constraints were met and thus not violated. The optimization algorithm has succeeded to meet the slip angle constraint. The optimization control still shows a good performance in controlling the vehicle to follow the optimum path. And shown in Fig. \ref{e4}(b) slip angle data for the rear tires at each time instant. It is shown that the slip angle is exactly between the maximum and minimum limits. The optimization algorithm has succeeded to meet the set constraint for the rear slip angles.

The simulation result and data previously discussed shows that the proposed optimization algorithm could successfully control the vehicle to maintain the ideal road and avoid an obstacle without violating any of the presented constraints for the vehicle motion. The optimization algorithm could successfully pass the steering angle constraint as shown in Fig. \ref{e3}(a), and braking ratio constraint as shown in Fig. \ref{e3}(b). The optimization control was also be able to control the vehicle to follow the ideal road and avoid an obstacle while minimizing the cost function to an optimum value as shown in Fig. \ref{e1}(b) and Fig. \ref{e2}(a). The proposed optimization algorithm has succeeded to meet the constraint for the front tires slip angle and to meet the constraint of the rear tires slip angle as shown in Fig. \ref{e4}. A summary of the evaluation done for the scenario test is shown in Table \ref{s8}.


\begin{table}[h!]

\def\arraystretch{1.5}
\caption{Scenario evaluation for avoiding an obstacle at a velocity of 10 m/s}
\centering
\begin{tabular}{||p{2cm}|p{2cm}|p{2cm}||}

 \hline
 Test & Value &Evaluation\\
 \hline
 $\text{Vehicle Motion}$ & - & $\text{Success}$\\
 $\text{Steering Angle}$ & $\in [-0.042, 0.06]$ & $\text{Success}$\\
 $\text{Braking Ratio}$ & $\in [-1, 1]$ & $\text{Success}$\\
 $\text{Front Slip Angle}$ & $\in [-0.6, 0.6]$ & $\text{Success}$\\
 $\text{Rear Slip Angle}$& $\in [-0.05, 0.1]$ & $\text{Success}$ \\
 \hline

\end{tabular}

\label{s8}
\end{table}	
During these tests the algorithm had achieved the goal by controlling the vehicle exactly at the ideal road without violating the lane safe margins. In some cases, the slip angle constraints were violated which could risk a slippage of the vehicle. Hardware applications can further control the slip angle to a limit where the vehicle is maintained on the ground. In other cases, an obstacle was presented to test the ability of the algorithm to avoid the obstacle while keeping the lane. The proposed algorithm achieved the goal and could avoid a smooth obstacle and rough obstacle. Another constraint that was implemented but not discussed is the maximum vehicle velocity. A constraint for the maximum velocity was set so that the optimization algorithm should not violate. In all the test cases, the algorithm was successful at maintaining the presented maximum velocity by manipulating the braking ratio of the vehicle. The braking ratio and steering angle showed a rapid behavior in some cases in order for the optimization to achieve the required goals. This rapid change could affect the comfort of the passenger and it could be further improved by applying it smoothly. More scenarios were tested but were not included in the work due to allowed space and limitations, therefore only interesting results were shown and discussed.

When compared with the work in \cite{20}, it was noticed that the proposed approach in this work shows more significant results regarding lane-keeping and collision avoidance systems’ performance. The control algorithm presented here could control a vehicle at a considerably high velocity through a mathematical modeling that wasn’t achievable in previous works. In addition to controlling a vehicle at high velocity, the proposed control algorithm was able to avoid a smooth obstacle, a steep obstacle, and keeping the vehicle in its road lane of different geometries including parabolic, inclined, and straight. Through analyzing mathematical performance of the presented work, it was found that all mathematical constraints, which were faced during the control and optimization process, were successfully met.

The previous work done in the field of automotive for passive safety systems was not found reliable in predicting the likelihood of having an accident \cite{ref1}. Another systems could predict the likelihood of having an accident but it was not integrated with a lane - keeping feature \cite{ref2}. Other lane - keeping systems were not integrated with collision avoidance feature \cite{20}. The proposed system could effectively avoid a collision based on the information given by a sensor and also keep the vehicle in lane if there was no information about an obstacle. The system could be implemented into action as a standalone decision making system that could take actions and override the driver commands.

\section{Conclusion and Future Prospects}

The proposed work has presented a full control on vehicle for lane - keeping and collision avoidance at low velocity and high velocity. The simulation results and evaluation demonstrated the capability of the proposed algorithm of controlling the vehicle at various scenarios and was successful. Our computations also showed that the vehicle could avoid a collision that was modeled as a smooth obstacle and rough obstacle. The results demonstrated that the optimization and control algorithm could also control the vehicle to stay in lane for straight, parabolic, and inclined lane models at low velocity. In addition, the proposed work showed the ability to control the vehicle at high velocities for inclined lane. 

The future work will concentrate on improving the computation speed to reach the level of a ”in the loop simulation”. Another aspect to improve the computation is to examine the influence of different tire models on the stability of the results. To reach the status of a in loop computation, a hardware implementation should be examined as a prototype. An open problem is under which parametric values of the obstacle a stable solution for the vehicle can be found. This kind of question needs to examine the influence of the structure of the obstacle on the optimization processes. Another aspect to consider for the future work is the effect of rapidly changing steering angle and braking ratio on passenger’s comfort and the trade-offs of the process. Some of these questions will be part of a future paper.

\section*{Acknowledgment}

This work was performed on the computational resource bwUniCluster funded by the Ministry of Science, Research and the Arts Baden-W{\"u}rttemberg and the Universities of the State of Baden-W{\"u}rttemberg, Germany, within the framework program bwHPC.

\ifCLASSOPTIONcaptionsoff
  \newpage
\fi



\bibliographystyle{IEEEtr}
%

%

\begin{biography}[{\includegraphics[width=1.2in,height=1.2in,clip,keepaspectratio]{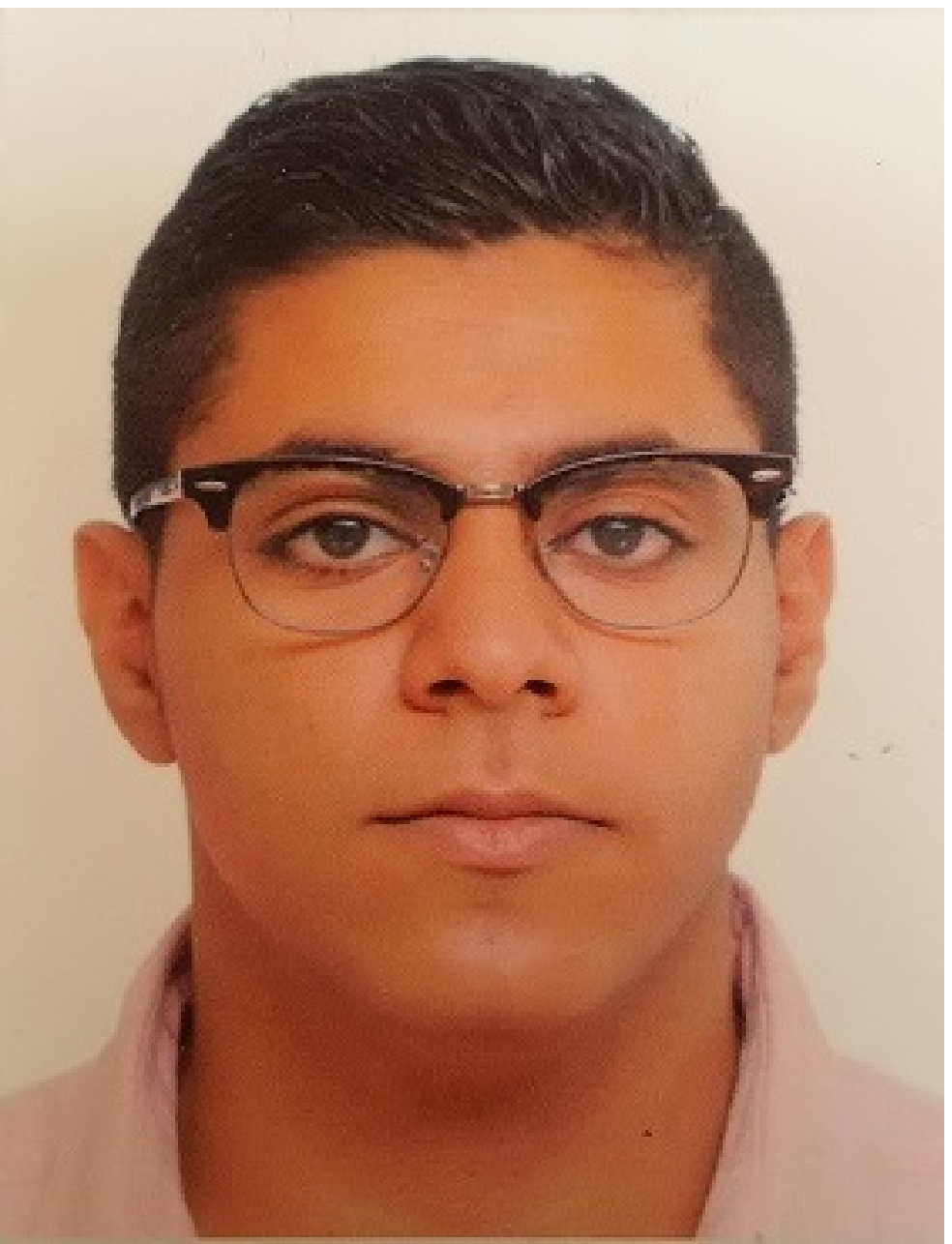}}]{Hazem M. Fahmy} received his BSc. degree in Electronics engineering from German University in Cairo, Egypt, in 2014 on Renewable Energy and Concentrated Solar Power. He received his MSc. degree in Electronics engineering from German University in Cairo, Egypt, in 2015 on Automotive Mathematical modeling. His research interest is in Automotive Applications, Embedded Systems, Renewable Energy and Mathematical modeling.
\end{biography}

\begin{biography}[{\includegraphics[width=4in,height=2in,clip,keepaspectratio]{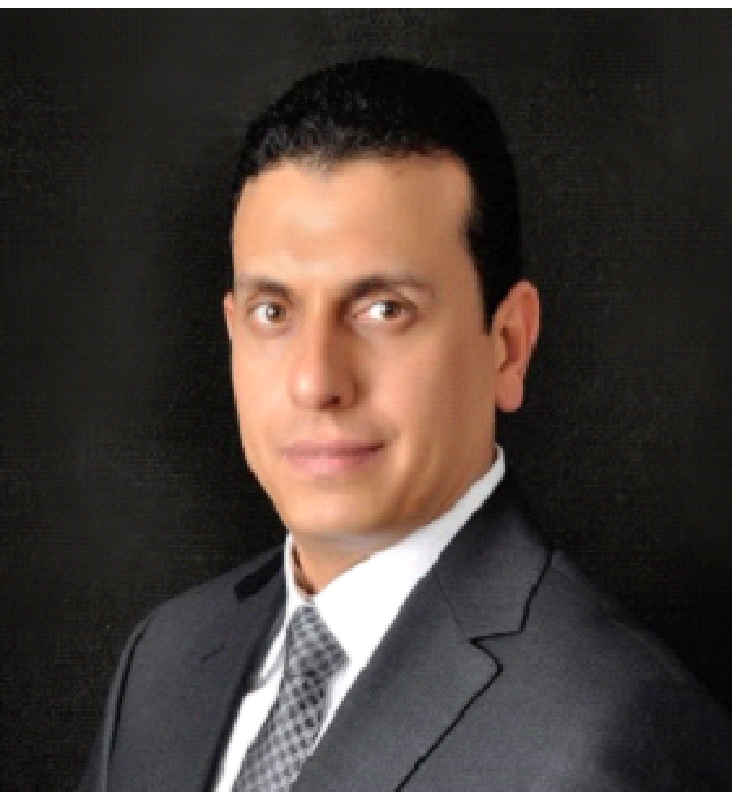}}]{Mohamed A. Abd El Ghany}  received the B.S. degree in electronics and communications engineering (with honors) and the Masters degrees in electronics engineering from Cairo University, Cairo, Egypt, in 2000 and 2006, respectively, and Ph.D degree in the area of high-performance VLSI/IC design from the German University, Cairo, Egypt in 2010. From 2003 to 2006, he was in National Space Agency of Ukraine, EGYPTSAT-1 project. From 2008 to 2009, he was an International Scholar at the Ohio State University, Electrical Engineering Dept., Columbus, USA. From 2012 to 2014, He awarded the Alexander von Humboldt Foundation Postdoctoral Fellowship, TU, Darmstadt, Germany. He is currently working as an Assistant Professor in German University in Cairo, Egypt. He is a project manager for two inernational projects between TU Darmstadt, Ruhr- Universitt Bochum and German University in Cairo. His research interest is in Network on Chip design and related circuit level issues in high performance VLSI circuits, clock distribution network design, digital application-specified integrated circuit design, SoC/NoC deisgn and verfication, low-power design and embedded system design. He is the author of about 30 papers, two book chapters, two book in the fields of high throughput and low-power VLSI/IC design and network on chip (NoC)/system on chip (SoC). He is a reviewer and program committee member of many IEEE international journals and conferences.
\end{biography}

\begin{biography}[{\includegraphics[width=3.2in,height=1.2in,clip,keepaspectratio]{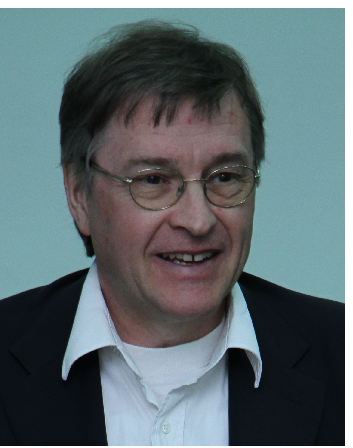}}]{Gerd Baumann} received his Dipl.-Phys. degree with honors in 1985, his PhD degree he earned with honors in 1988 in the field of Mathematical
Physics from the University of Ulm, Ulm, Germany, respectively. In 1993 he received his Venia Legendi for Theoretical Physics and was first appointed as Professor in 1998 at the University of Ulm, Ulm, Germany. He received different awards from the Hans Voith Foundation, an award of ”The Baden- Wuerttemberg Employers’ Association of the Metal Industry e. V.”, Germany, an award of the University Association Ulm, a Scholar Grant by Wolfram Research, US, and the Merckle Research Award, Germany. In 1989 he was visiting as a researcher the Center for Nonlinear Studiesin at the Los Alamos National Laboratory, US. Starting in 2000 he was working beside his University activities as Head of Advanced Engineering and Consultant with major international companies in Automotive, Medicine, Software Engineering, and Biology. In 2004 he was appointed Professor and Head of the Mathematics Department at the German University in Cairo, Cairo, Egypt. He is currently managing a bilateral international research cooperation funded by BMBF between Germany and Egypt in the field of lightweight and degradable materials in biological applications. His research interests are in applied hybrid Mathematics for industrial applications. He is author of numerous papers and books aiming on practical applications in industrial environments. As a reviewer and committee member at DAAD and in different international journals he is actively participating in new developments in science and engineering. He is a member of SIAM and AMS, US.
\end{biography}





\end{document}